\pdfoutput=1
\RequirePackage{ifpdf}
\ifpdf 
\documentclass[pdftex]{sigma}
\else
\documentclass{sigma}
\fi

\usepackage[all]{xy}
\usepackage{tikz-cd}
\usepackage{stmaryrd}
\usepackage{mathtools}

\numberwithin{equation}{section}

\newtheorem{Theorem}{Theorem}[section]
\newtheorem*{Theorem*}{Theorem}
\newtheorem{Corollary}[Theorem]{Corollary}
\newtheorem{Lemma}[Theorem]{Lemma}
\newtheorem{Proposition}[Theorem]{Proposition}
 { \theoremstyle{definition}
\newtheorem{Definition}[Theorem]{Definition}

\newtheorem{Example}[Theorem]{Example}
\newtheorem{Examples}[Theorem]{Examples}
\newtheorem{Remark}[Theorem]{Remark} }

\newtheorem{innercustomthm}{Theorem}
\newenvironment{customthm}[1]
{\renewcommand\theinnercustomthm{#1}\innercustomthm}
{\endinnercustomthm}

\newcommand{\R}{{\mathbb{R}}}
\newcommand{\C}{{\mathbb{C}}}

\newcommand{\g}{\mathfrak{g}}

\newcommand{\X}{\mathfrak{X}}
\newcommand{\ra}{\rangle}

\newcommand{\CMC}{C^{\infty}(M, \C)}
\newcommand{\Diff}{\mathrm{Diff}}

\newcommand{\T}{\mathbb{T}}

\newcommand{\cA}{\overline{A}}

\newcommand{\aut}{\mathfrak{aut}}

\DeclareMathOperator{\jc}{Jac}
\DeclareMathOperator{\type}{type}

\DeclareMathOperator{\order}{order}

\DeclareMathOperator{\Aut}{Aut}

\DeclareMathOperator{\im}{Im}
\DeclareMathOperator{\re}{Re}
\DeclareMathOperator{\rk}{rank}
\DeclareMathOperator{\cork}{corank}

\DeclareMathOperator{\pr}{pr}

\DeclareMathOperator{\real}{re}
\DeclareMathOperator{\rerk}{real-rank}
\DeclareMathOperator{\class}{class}

\begin{document}

\newcommand{\arXivNumber}{2401.05274}

\renewcommand{\PaperNumber}{044}

\FirstPageHeading

\ShortArticleName{On Complex Lie Algebroids with Constant Real Rank}

\ArticleName{On Complex Lie Algebroids with Constant Real Rank}

\Author{Dan AGUERO}

\AuthorNameForHeading{D.~Aguero}

\Address{Scuola Internazionale Superiore di Studi Avanzati - SISSA,\\
Via Bonomea, 265, 34136 Trieste, Italy}
\Email{\href{mailto:dagueroc@sissa.it}{dagueroc@sissa.it}}

\ArticleDates{Received September 30, 2024, in final form June 02, 2025; Published online June 13, 2025}

\Abstract{We associate a real distribution to any complex Lie algebroid that we call {\em distribution of real elements} and a new invariant that we call {\em real rank}, given by the pointwise rank of this distribution. When the real rank is constant, we obtain a real Lie algebroid inside the original complex Lie algebroid. Under another regularity condition, we associate a complex Lie subalgebroid that we call {\em the minimal complex subalgebroid}. We also provide a local splitting for complex Lie algebroids with constant real rank. In the last part, we introduce the {\em complex matched pair} of skew-algebroids; these pairs produce complex Lie algebroid structures on the complexification of a vector bundle. We use this operation to characterize all the complex Lie algebroid structures on the complexification of real vector bundles.}


\Classification{53D20}

\section{Introduction}
Lie algebroids were introduced by Pradines \cite{pradines3lie}, and nowadays have become a useful tool for many areas of differential geometry and mathematical physics. This influenced their study and made a fast development of the area. However, the study of Lie algebroids can be divided into real, holomorphic, algebraic and complex Lie algebroids, with real Lie algebroids the most studied. Complex Lie algebroids gained more attention after the introduction of generalized complex structures \cite{gualtieri:2007, hitchin:2004} and now we know that many geometrical structures are particular examples of complex Lie algebroids, for example, generalized complex structures, complex Dirac structures \cite{aguero2022complex}, exceptional complex structures \cite{tennyson2021exceptional}, $B_n$-generalized complex structures \cite{rubio2014generalized}, etc.

The theory of real and holomorphic Lie algebroids is currently more developed than the theory of complex Lie algebroids. It is known their local structure, the obstructions for its integrability (Lie~III) \cite{crainic2003integrability,laurent2009integration}, etc. Despite the general theory of generalized complex structures and complex Dirac structures being known \cite{aguero2022complex, bailey2013local}, it is still missing a general description of complex Lie algebroids including its local description. About the integrability of complex Lie algebroids there is not much known, some ideas can be found in \cite{weinstein2007integration}.

In this article, we study complex Lie algebroids, by introducing a real distribution that we call {\em distribution of real elements} and an invariant that we call the {\em real rank} which is just the rank of the aforementioned distribution. When the real rank is constant the distribution of real element is a real Lie algebroid and so we can associate a real foliation to the complex Lie algebroid.
In~the~philosophy of associating real geometrical structures, the Lie algebroid of real elements is an analog to the Poisson and Dirac structures associated to generalized complex and complex Dirac structures, respectively. We associate two other invariants: the type and the class; and introduce the minimal complex Lie subalgebroids and the family of minimal complex Lie algebroids. Our main examples are involutive structures, complexifications and complex Poisson bivectors. We provide a local splitting theorem for complex Lie algebroids with constant real rank. We also study complex Lie algebroid structures on complexifications of real vector bundles, we introduce the complex matched pairs and an operation associated with them that we call complex sum. We prove that there is a one-to-one correspondence between complex Lie algebroid structures on complexifications and complex matched pairs.

The content of this article is organized as follows. In Section~\ref{CAVB}, we recall some properties of the symmetries of complex vector bundles and complex anchored vector bundles. Since the author has not found the proofs of some of the statements, they are provided. After that, we focus on complex anchored vector bundles (CAVBs). We introduce the {\em distribution of real elements} and its rank which we call {\em real rank}. When the distribution of real elements has constant rank it is a real anchored vector bundle, in fact, it is involutive in case the CAVB is involutive. In~Section~\ref{CLA}, we study complex Lie algebroids (CLAs), under the assumption of constant real rank. We introduce {\em the Lie algebroids of real elements} denoted by $(A^{\real},[\cdot,\cdot]^{\real},\rho^{\real})$ and we study some of its properties. This real Lie algebroid allows us to associate a real foliation to a CLA with constant real rank. Since we are dealing with complex distributions, we can have many types of regularity rather than the usual regularity (the image of the anchor map is a~regular distribution), so we introduce strongly regular CLAs. We also study the distribution given by $A_{\min}=A^{\real}+{\rm i}A^{\real}$ and whenever $A_{\min}$ is regular it is automatically a CLA, which we call {\em the~minimal complex Lie subalgebroid}. We introduce a new class of CLA that we call {\em minimal CLA}, which are CLAs that coincide with their minimal complex subalgebroid. We also introduce an $\mathbb{N}$-valued function called the {\em class} that takes a point and returns the dimension of the leaf passing through this specific point. CLAs of class $0$ at a given point play an important role in the splitting theorem for CLAS. In the last part of the section, we study the conjugate of a~CLA.

In Section~\ref{locstr}, we study the local structure of CAVBs and CLAs with constant real rank using the techniques of \cite{bursztyn2019splitting}. We provide a proof of the splitting theorem for involutive CAVBs with constant real rank and we end the section with the local description of CLAs with constant real rank, obtaining the following:
\begin{customthm}{4.7}

  Let $(A,[\cdot,\cdot],\rho)$ be a CLA, $m\in M$ and $\Delta=\rho(A)\cap TM$. Assume that $(A,[\cdot,\cdot],\rho)$ has constant real rank in a neighborhood of $m$. Consider a submanifold \smash{$N\xhookrightarrow{\iota}M$} such that $\Delta_{m}\oplus T_m N=T_m M$ and denote by $P=\Delta_{m}$. Then, there exists a neighborhood $U$ of $m$ such that $A|_{U}$ is isomorphic to $ A'\times T_{\C}P$, where $(A',[\cdot,\cdot]',\rho')$ is a CLA of class $0$ at $m$.
\end{customthm}
The property of $A'$ having class $0$ at $m$ already appears in the local description of generalized complex structures and complex Dirac structures.

As a main application, we recover the Frobenius--Nirenberg theorem \cite[Theorem 1]{nirenberg1958complex}.

In Section~\ref{cxsum}, we introduce the {\em complex matched pairs of skew-algebroids}. Given two skew-algebroids $(A, [\cdot,\cdot]_1, \rho_1)$ and $(A, [\cdot,\cdot]_2, \rho_2)$ over the same vector bundle $A\to M$, we say that they form a complex matched pair if the brackets are compatible and the Jacobiator of the brackets is the same. We observe that these pairs produce a complex Lie algebroid structure on the complexification of a real vector bundle which we call {\em complex sum}. As an application, we prove the following characterization:
 \begin{customthm}{5.6}

 Let $A$ be a real vector bundle over $M$. Any structure of CLA on $A_{\C}$ comes from a complex sum of a complex matched pair of skew-algebroids\footnote{Briefly speaking, a skew-algebroid is a Lie algebroid with a bracket that does not necessarily satisfy the Jacobi identity, see Appendix \ref{preliminaries0}.} over $A$.
\end{customthm}
We also study the CLA associated to complex vector fields and describe its distribution of real elements. In this case, we see that the assumption of constant real rank is inefficient for the local description of general complex vector fields. However, allowing a ``singular associated real Lie algebroid'' could provide a solution to this problem, this is the subject of future work. In the end of the section, we study the CLA associated to a complex Poisson bivector from the perspective of the complex sum.
In the end of this article, we present two appendices: Appendix \ref{preliminaries0}, where we recall some results and properties of real anchored vector bundles, vertical lifts and real Lie algebroids, and Appendix \ref{preliminarieslocstr}, where we recall some results and properties about Euler-like vector fields and normal vector bundles.
\subsection*{Notation and convention}
In this article, we only work with smooth manifolds. Given a vector bundle $A$ over a smooth manifold $M$, a {\em distribution} is an assignment $m\in M\mapsto R_m\subseteq A_m$, where $R_m$ is a vector subspace of the fiber $A_m$; when any $v\in R_m$ can be extended to a local section of $A$ taking values in $R$, we say that the distribution is {\em smooth}. Each distribution has associated a $\mathbb{N}$-valued function given by the assignment $\rk\colon m\in M\mapsto\dim R_m$, which is called {\em rank}. A distribution is said to be {\em regular} if it is a vector bundle.

Given a distribution $E\subset A_\C$, where $A_{\C}$ is the complexification of a real vector bundle $A\to M$, we denote the space of real elements of $E$ by $\re E:=E\cap A$ and call them the {\em real part of $E$}. Given a smooth map $\varphi\colon M\to N$, we denote the complexification of $T\varphi$ by $T_{\C}\varphi$.

Throughout this article, a real anchored vector bundle is referred to by RAVB and a real Lie algebroid by RLA. In Appendix~\ref{preliminaries0}, we recall some properties of RAVBs and RLAs.

Let $N$ be a submanifold of $M$, the normal bundle of $N$ is denoted by
\[
\nu(M,N)=TM|_{N}/TN,
\]
when the context is clear we denote it by $\nu_N$. In Appendix \ref{preliminarieslocstr}, we recall some properties of normal vector bundles.

\section{Complex anchored vector bundles}\label{CAVB}
\subsection{Complex anchored vector bundles}
A complex vector bundle $E\to M$ is equipped with a canonical bundle map $j_E\in \Diff(E)$ defined by $j_E(e)={\rm i}\cdot e$, the multiplication by the imaginary number ${\rm i}$. The {\em group of $($complex$)$ automorphisms of $E$} is
 \[
 \Aut(E)=\{(\Phi,\varphi)\in \Aut(E_{\R})\mid \Phi\circ j_E=j_E \circ \Phi\ \text{and}\ \Phi\ \text{covers}\ \varphi\in \Diff(M)\},
 \]
 where $E_{\R}$ denotes the realification of $E$, and $\Aut(E_{\R})$ is the automorphism group of $E_{\R}$ as a~real vector bundle.
The map $j_E$ acts on $\Aut(E)$ in the natural way. Consider a 1-parameter subgroup $\{\Phi_t\}_{t\in \R}$ of $\Aut(E)$. Since $\Phi_t\in \Aut(E_{\R})$ and $\Phi_t\circ j_E=j_E\circ \Phi_t$, the vector field \smash{$\hat{X}=\frac{{\rm d}}{{\rm d}t}\Phi_t|_{t=0}\in \mathfrak{aut}(E_{\R})$} satisfies \smash{$(j_E)_* \hat{X}=\hat{X}$}. So, the {\em Lie algebra of $($complex$)$ infinitesimal automorphisms of $E$} is given by
 \[\mathfrak{aut}(E)=\bigl\{\hat{X}\in \mathfrak{aut}(E_{\R})\mid (j_E)_* \hat{X}=\hat{X}\bigr\}\subseteq \mathfrak{aut}(E_{\R})\subseteq \mathfrak{X}(E).\]
 Another characterization of $\aut(E)$ is the following.

 \begin{Lemma} The Lie algebra $\mathfrak{aut}(E)$ is given by the pairs $(L,X)$ such that $X\in \mathfrak{X}(M)$ and $L\colon \Gamma (E)\to \Gamma (E)$ is a $\C$-linear derivation with respect to $X$.
 \end{Lemma}
\begin{proof}
Since $\widetilde{X}\in \aut(E)\subseteq \aut (E_{\R})$, it has associated a derivation given by the pair $(L, X)$, where $L\colon \Gamma(E_{\R})\to \Gamma(E_{\R})$ is $\R$-linear and \smash{$X=\pr_* \widetilde{X}\in \X(M)$}, defined as \smash{$L(\sigma)^{\uparrow}=[\sigma^{\uparrow}, \widetilde{X}]$}, ${\forall \sigma\in \Gamma(E)}$, where $\uparrow$ denotes the vertical lift, see Appendix \ref{preliminaries0}. Note that $L$ is $\C$-linear if and only if $L({\rm i}\sigma)={\rm i}L(\sigma)$ and thus
\begin{align*}
L({\rm i}\sigma)^{\uparrow}&{}=\bigl[({\rm i}\sigma)^{\uparrow}, \widetilde{X}\bigr]=\bigl[(j_{E})_* \sigma^{\uparrow}, \widetilde{X}\bigr]=\bigl[(j_{E})_* \sigma^{\uparrow},(j_{E})_* \widetilde{X}\bigr]\\
&{}=(j_E)_* \bigl[\sigma^{\uparrow}, \widetilde{X}\bigr]=(j_{E})_* \bigl(D(\sigma)^{\uparrow}\bigr)=({\rm i}L(\sigma))^{\uparrow}.
\end{align*}
And by injectivity of the ``$\uparrow$''-map, we have that $L({\rm i}\sigma)={\rm i}L(\sigma)$ and the lemma holds.
\end{proof}

\begin{Definition}
 A {\em complex anchored vector bundle (CAVB)} is a pair $(A,\rho)$, where $A$ is a complex vector bundle over $M$, and $\rho\colon A\to T_{\C}M$ is a $\C$-linear bundle map. Let $(A,\rho_A)$ and $(B,\rho_B)$ be two CAVBs over $M$ and $N$, respectively. A {\em morphism of complex anchored vector bundles} consists of a bundle map $\Phi\colon A\to B$ and a map $\varphi\colon M\to N$ such that the following diagram commutes:
\[\xymatrix{A\ar[d]^{\rho_A}\ar[r]^{\Phi} & B\ar[d]^{\rho_B}\\ T_{\C}M\ar[r]^{T_{\C}\varphi} & T_{\C}N.
}\]

\end{Definition}
Given a CAVB $(A, \rho)$, the anchor map decomposes $\rho=\rho_1+{\rm i}\rho_2$. So we obtain two real anchored vector bundles $(A_{\R},\rho_1)$ and $(A_{\R}, \rho_2)$. By the $\C$-linearity of $\rho$, we have that
\[\rho_1 \circ j_A=-\rho_2.\]
\begin{Proposition}
The map $j_A$ is an isomorphism of RAVBs between $(A, \rho_1)$ and $(A,-\rho_2)$.
\end{Proposition}
 The {\em group of automorphisms of a CAVB} $(A,\rho)$ is given by
 \[\Aut(A,\rho)=\{(\Phi, \varphi)\in \Aut(A)\mid \rho\circ\Phi=T_{\C}\varphi\circ \rho\}.\]
 Note that $\rho\circ\Phi=T_{\C}\varphi\circ \rho$ if and only if $\rho_1\circ \Phi=T\varphi\circ \rho_1$ and $\rho_2\circ \Phi=T\varphi\circ \rho_2$. So
 \[\Aut(A,\rho)=\Aut(A_{\R},\rho_1)\cap \Aut(A_{\R},\rho_2)\cap \Aut (A).\]
 The {\em $($complex$)$ infinitesimal automorphisms of $(A,\rho)$} are denoted by $\aut(A,\rho)$.
\begin{Definition}\label{cx_tangent_extension}
   The {\em complex tangent extension} of a vector field $X\in \X(M)$, is the vector field $X_{T_{\C}}\in \mathfrak{X}(T_{\C}M)$ defined by the one-parameter subgroup \smash{$T_{\C}\varphi_t ^{X}$}, where \smash{$\varphi_t ^{X}$} is the local flow of $X$.
\end{Definition}

For a CAVB $(A,\rho)$, we have that
\begin{align*}
\mathfrak{aut}(A,\rho)&{}=\aut(A_{\R},\rho_1)\cap \aut(A_{\R},\rho_2)\cap \aut(A)\\
&{}=\bigl\{\widetilde{X}\in \X(A)\mid (j_E)_* \hat{X}=\hat{X},\, \widetilde{X}\sim_{\rho} X_{T_{\C}},\, \text{where}\ \widetilde{X}\sim_{\pr} X\bigr\}.
\end{align*}

In terms of the differential operators: a pair $(L,X)$, where $L\colon \Gamma(A)\to \Gamma(A)$ is a $\C$-linear differential operator and $X\in \mathfrak{X}(M)$, belongs to $\mathfrak{aut}(A,\rho)$ if and only if $\rho(L(\tau))=[X,\rho(\tau)]$
or equivalently
\begin{equation*}
\rho_i(L(\tau))=[X, \rho_i(\tau)]\qquad\text{for $i=1,2$}.
\end{equation*}
\begin{Definition}
A real or complex anchored vector bundle $(A,\rho)$ is {\em involutive} if $\rho(\Gamma(A))$ is a~Lie subalgebra of $\Gamma(TM)$ or $\Gamma(T_{\C}M)$, respectively.
\end{Definition}

Consider the decomposition of the anchor map $\rho$ in its real and imaginary parts ${\rho=\rho_1+{\rm i}\rho_2}$. Note that the involutivity of $(A,\rho)$ is not directly related to the involutivity of $(A_{\R}, \rho_1)$ and $(A_{\R}, \rho_2)$. Even if both $(A_{\R}, \rho_1)$ and $(A_{\R}, \rho_2)$ are involutive, we cannot ensure the involutivity of~$(A,\rho)$ and vice versa. To solve this issue, we shall introduce the distribution of real elements in the next subsection.

\subsection{Distributions associated to CAVBs}
A CAVB $(A,\rho)$ has associated two real distributions~in~$TM$:
 \[\Delta=\re \bigl(\rho(A)\cap \overline{\rho(A)}\bigr)\qquad\text{and} \qquad D=\re \bigl(\rho(A) + \overline{\rho(A)}\bigr).\]
Observe that $D$ is always a smooth distribution, while $\Delta$ is not necessarily smooth. However, we have the following.

 \begin{Lemma}\label{imageanchor}
 If $D$ is a regular distribution, then $\Delta$ is a smooth distribution.
 \end{Lemma}
\begin{proof}
 Consider the map $\rho_2\colon A\to TM$, which is a real bundle map. Since the image of $\rho_2$ is~$D$, its kernel $\ker \rho_2$ is a vector bundle. Thus, $\Delta=\rho(\ker \rho_2)$ is smooth.
\end{proof}

\begin{Definition}
Given a CAVB $(A,\rho)$ define the {\em distribution of real elements of $(A,\rho)$} as
\[A^{\real} =\{e\in A\mid \rho_2(e)=0\},\]
which we called {\em distribution of real elements of $(A,\rho)$}. In general, $A^{\real}$ is a distribution in $A$ (not necessarily smooth). When $A^{\real}$ is regular, the pair $(A^{\real},\rho^{\real})$ is a RAVB, where $\rho^{\real}=\rho_1|_{A^{\real}}$. The $\mathbb{N}$-valued function 
\[\operatorname{\text{real-rank}} A\colon \ m\in M\mapsto \dim_{\R} A^{\real}|_m\]
is called {\em real rank of~$A$}.
\end{Definition}
Constant real rank is equivalent to $A^{\real}$ being a vector bundle and by the proof of Lemma \ref{imageanchor}, it is equivalent to the regularity of $D$ as a distribution. Indeed, we can see that
\begin{equation}\label{sum_identity}
  {\operatorname{\text{real-rank}} A}+\rk_{\R} D=\rk_{\R} A_{\R}.
\end{equation}
\begin{Lemma}\label{restriction_inf_aut}
   Assume that $(A,\rho)$ has constant real rank. If $(L,X)\in \aut(A,\rho)$, then $L(\Gamma(A^{\real}))\subseteq\Gamma(A^{\real})$ and $(L|_{A^{\real}}, X)\in \aut(A^{\real},\rho^{\real})$.
\end{Lemma}
\begin{proof}
  Straightforward.
\end{proof}

 \begin{Lemma}\label{cavb_inv}
 If $(A,\rho)$ is an involutive CAVB with constant real rank, then $(A^{\real}, \rho^{\real})$ is involutive.
 \end{Lemma}
\begin{proof}
Let $\alpha, \beta \in \Gamma (A^{\real})$, equivalently, $\rho (\alpha)$ and $\rho (\beta)$ are real, and so $[\rho (\alpha), \rho (\beta)]$. By involutivity of $A$, there exists $\gamma \in \Gamma (A)$ such that $[\rho (\alpha), \rho (\beta)]=\rho (\gamma)$. So $\rho (\gamma)$ is real and consequently~${\gamma\in \Gamma(A^{\real})}$.
\end{proof}

The image of the anchor map of an involutive RAVB is integrable, see for example \cite{bursztyn2019splitting}. In~the~case of a CAVB $(A,\rho)$, it follows from the proof of Lemma \ref{imageanchor} and by Lemma \ref{cavb_inv} that the real part of $\rho(A)$ is integrable.

The pullback of complex anchored vector bundle $(A, \rho) $ over $M$ along the map $\varphi\colon N\to M$ is defined in the same way as the real case, as the fibered product:
\[
\xymatrix{\varphi^! A\ar[d]\ar[r] & A \ar[d]^{\rho}\\
T_{\C}N\ar[r]^{T_{\C}\varphi} & T_{\C}M.}
\]

Note that
\begin{align*}
\varphi^! A&{}=\bigl\{(a, X_1+{\rm i}X_2) \mid  (a,X_1)\in \varphi^! A_1\ \text{and}\ (a,X_2)\in \varphi^! A_2\bigr\}\\
&{}=\varphi^! A_1 \times_{\pr_A} {\rm i}\varphi^! A_2\subseteq A\times T_{\C}N,
\end{align*}
where $A_1=(A_{\R}, \rho_1)$ and $A_2=(A_{\R}, \rho_2)$, and $\pr_A\colon A\times T_{\C}N\to A$ is the obvious projection.
As~a~consequence, we can see that
\begin{enumerate}\itemsep=0pt
\item $\varphi^!T_{\C}M=T_{\C}N$.
\item If $N=M\times Q$, then $\pr_M ^! A=A\times T_{\C}Q$.
\item If \smash{$N\xhookrightarrow{\iota} M$} is a submanifold transversal to $\rho$, then $N$ is also transversal to $\rho_1$ and $\rho_2$. Moreover,
\[
\bigl(\iota^! A\bigr)_{\R}=\rho^{-1}(T_{\C}N)_{\R}=\{a\in A_{\R}\mid \rho(a)=\rho_1(a)+{\rm i}\rho_2(a)\in T_{\C}N\}=\iota^{!}A_1\cap \iota^{!}A_2.
\]
\item If $\rho$ is injective, then $\varphi^! A=(T_{\C}\varphi)^{-1}(A)$.
\end{enumerate}

 The following lemma is straightforward.

\begin{Lemma}\label{repartbackward}
Let $(A,\rho)$ be a CAVB and a map $\varphi\colon N\to M$. Then,
\begin{align*} 
&\bigl(\varphi^! A\bigr)^{\real}=\varphi^! (A^{\real})\times_{\pr_A}{\rm i}(A\times \ker T\varphi),\\
&\Delta_{\varphi! A}=\im(T\varphi)\cap \Delta_{A}.
\end{align*}
\end{Lemma}

A submanifold \smash{$N\xhookrightarrow{\iota} M$} is called {\em transversal to the CAVB} if
\[\rho(A)+T_{\C}N=T_{\C}M\]
(or the stronger condition $\Delta+TN=TM$). By transversality, \smash{$\iota^! A$} is smooth and so is a CAVB. Note that transversality with respect to $\rho$ implies transversality with respect to $A_1$ and $A_2$.

A direct consequence of the Lemma \ref{repartbackward} is that \smash{$\bigl(\iota^! A\bigr)^{\real}=\iota^! (A^{\real})$} for submanifolds \smash{$N\xhookrightarrow{\iota} M$} (not necessarily transversal to the anchor map).
The following lemma appeared originally in~\cite{bursztyn2019splitting}, and the proof is similar in the complex case.
\begin{Lemma}
   Suppose $(A,\rho)$ is an involutive CAVB over $M$ and $\varphi\colon N\to M$ a map transversal to $\rho$. Then $\varphi^! A$ is an involutive CAVB.
\end{Lemma}

 \section{Complex Lie algebroids}\label{CLA}
\begin{Definition}
A {\em complex Lie algebroid $($CLA$)$} over a manifold $M$ is a complex anchored vector bundle $(A, \rho)$ over $M$ together with a Lie bracket on $\Gamma(A)$ satisfying the Jacobi identity
\[[\cdot,\cdot]\colon \ \Gamma(A)\times \Gamma(A)\to\Gamma(A)\]
 and with a $\C$-bundle map $\rho\colon A\to T_{\C}M$, that we call the anchor map, preserving brackets and satisfying the Leibniz identity
 \[[\alpha,f \beta]=f[\alpha,\beta]+\rho(\alpha)(f)\beta\]
 for all $\alpha,\beta\in \Gamma(A)$ and for any function $f\in C^{\infty}(M,\C)$. A {\em complex skew-algebroid} is when the Jacobi identity is not necessarily satisfied, and an {\em almost complex Lie algebroid $($almost CLA$)$} is a complex skew-algebroid where the anchor preserves the brackets.
\end{Definition}

\begin{Remark}
Observe that CLAs are defined on the category of smooth real manifolds and not on the category of holomorphic manifolds.
\end{Remark}
Morphisms between CLAs are defined in the exact same way as RLAs, and so are their automorphisms.
\begin{Examples}\label{examples_0}\normalfont Some examples of CLAs are the following:\samepage
\begin{enumerate}\itemsep=0pt
\item An {\em involutive structure} $E$ is an involutive vector subbundle of $T_{\C} M$, see \cite{treves2014hypo}. The triple $(E, [\cdot,\cdot], \rho)$, where $[\cdot,\cdot]$ is the complexification of the Lie bracket of vector fields and $\rho$ is the inclusion map of $E$ in $T_{\C}M$, defines a CLA. In particular: complex structures, transverse holomorphic structures and CR structures define CLAs.

We recall that a CR structure is a pair $(I,H)$, where $H$ is a subbundle of $TM$ and ${I\colon H\to H}$ is a bundle map satisfying that $I^2=-{\rm Id}$ and the following: if ${X, Y\in H}$,~then
\[
  [IX,Y]+[X, IY]\in H\qquad \text{and}\qquad I([IX,Y]+[X, IY])=[IX,IY]-[X, Y].
\]
In this case, the triple $(E, [\cdot,\cdot], \iota)$ is a CLA, where $E=\ker(I_{\C}-{\rm i}{\rm Id})\subseteq T_{\C}M$, the bracket~$[\cdot,\cdot]$ is the complexification of the Lie bracket of vector fields and $\iota\colon E\to T_{\C}M$ is the inclusion map. The case of complex and transverse holomorphic structures is similar.

\item Complex Dirac structures \cite{aguero2022complex}, generalized complex structures \cite{gualtieri:2004}, $B_n$-generalized complex structures \cite{rubio2014generalized} and exceptional generalized complex structures \cite{tennyson2021exceptional}. Let $L$ be a complex Dirac structure, that is, a lagrangian subbundle $L$ of $(TM\oplus T^*M)_{\C}$ which is closed under the Courant--Dorfman bracket $\llbracket\cdot,\cdot\rrbracket$. The Courant--Dorfman bracket when restricted to the sections of an involutive isotropic subbundle become a Lie bracket. Thus, the CLA structure on $L$ is given by the triple $(L,\llbracket \cdot,\cdot\rrbracket|_{\Gamma(L)\times \Gamma(L)}, \pr_{T_{\C}M}$). It is similar for the other aforementioned structures.
\item A holomorphic Lie algebroid $(A, [\cdot,\cdot],\rho)$ has associated two CLAs, $(A, [\cdot,\cdot]_{1,0},\rho_{1,0})$ and $(A, [\cdot,\cdot]_{0,1},\rho_{0,1})$ defined in \cite[Section~3.4]{laurent2008holomorphic}.

\item A bundle of complex Lie algebras is a complex vector bundle $p\colon E\to M$ where each fiber~$E|_m$ is equipped with a complex Lie algebra bracket $[\cdot,\cdot]_m$ and satisfies the following: if~${s_1, s_2\in \Gamma(E)}$, then the assignment
\[
  m\in M\to [s_1,s_2]|_m:= [s_1|_m,s_2|_m]_m\in E|_m
\]
is a smooth section of $E$. Then, the triple defined by $(E,[\cdot,\cdot], \rho=0)$ is a CLA.
Conversely, any CLA $(A,[\cdot,\cdot],\rho)$ having $\rho=0$ is a bundle of complex Lie algebras (the proof is almost the same as in the real case, see for example \cite[Section~16.2]{da1999geometric}).
\item The complexification $(A_{\C}, [\cdot,\cdot]_{\C},\rho_{\C})$ of a RLA $(A, [\cdot,\cdot],\rho)$ defines a CLA.
\item Given a complex Lie algebra $\mathfrak{g}$ we call a Lie algebra homomorphism
$\mathcal{X}\colon \mathfrak{g}\to \Gamma(T_{\C}M)$ a~{\em complex infinitesimal action}. The vector bundle $L(\g)=M\times \g$ is a CLA with anchor map~$\mathcal{X}$ and bracket
\[[u,v](m)=[u(m), v(m)]+\bigl(L_{\mathcal{X}(u(m))}v\bigr)(m)-\bigl(L_{\mathcal{X}(v(m))}u\bigr)(m),\]
where $u$, $v$ and $[u,v]$ are sections of $L(\g)$ taken as maps from $M$ to $\g$.

\item
A complex Poisson bivector is a bivector $\pi \in \Gamma\bigl(\wedge^{2}T^*_{\C}M\bigr)$ satisfying $[\pi,\pi]=0$. The triple $\bigl(T^*_{\C}M, [\cdot, \cdot]_{\pi}, \pi\bigr)$ is a CLA, where the bracket $[\cdot, \cdot]_{\pi}$ is defined as
\begin{equation*}
[\alpha,\beta]_{\pi}=L_{\pi(\alpha)}\beta-L_{\pi(\beta)}\alpha- {\rm d}\pi(\alpha,\beta)
\end{equation*}
and anchor $\pi\colon T^*_{\C}M\to T_{\C}M$. We denote this complex Lie algebroid by \smash{$\bigl(T_{\C}^* M\bigr)_{\pi}$}. This is the complex parallel of the Lie algebroid associated to a Poisson structure.
\item A complex vector field $Z=X_1+{\rm i}X_2\in \Gamma(T_{\C}M)$ defines a structure of CLA over the bundle~${M\times \C}$ that we denote by $A_{Z}$. The anchor map is given by
\begin{align*}
 \rho\colon \ \Gamma(M\times\C)=C^{\infty}(M,\C)\to T_{\C}M,\qquad \rho(f)=fZ
\end{align*} and the bracket by
\[[f,g]=f L_Z (g)-g L_Z (f).\]
\end{enumerate}
\end{Examples}
\subsection{Associated RLA}
On one hand, a RLA have associated a distribution given by the image of the anchor map. On the other hand, the image of the anchor map of a CLA is a complex distribution of $T_{\C}M$. To~bypass this issue, we associate two real distributions:
\begin{equation}\label{real_distributions}
\Delta=\re \bigl(\rho(A)\cap \overline{\rho(A)}\bigr)\qquad\text{and} \qquad D=\re \bigl(\rho(A) + \overline{\rho(A)}\bigr).
\end{equation}
Note that $D$ is a smooth distribution, while $\Delta$ is not smooth in general, see \cite[Example 6.1]{aguero2022complex}.

Given a CLA $(A, [\cdot,\cdot], \rho)$, decompose $\rho=\rho_1+{\rm i}\rho_2$. Contrary to what happens with CAVBs, neither $\mathcal{A}_1=(A_{\R}, [\cdot,\cdot], \rho_1)$ nor $\mathcal{A}_2=(A_{\R}, [\cdot,\cdot],\rho_2)$ necessarily define RLAs, since $\rho_1$ and $\rho_2$ do not preserve the bracket (and so they do not satisfy the Leibniz identity).

Actually, the natural candidate for $\mathcal{A}_2$ is $(A_{\R}, {\rm i}[\cdot,\cdot], -\rho_2)$, here we write $-\rho_2$ in order to avoid~${\rm i}\rho_2$. The failure of the Leibniz identity for $\mathcal{A}_1$ and $\mathcal{A}_2$ is controlled in the following way:
  \begin{gather}\label{fail_leib_re}
  [a, fb]_1-(f[a,b]_1+(\rho_1(a)f)b)={\rm i}(\rho_2(a)f)b,
  \\ \label{fail_leib_im}
  [a, fb]_2-(f[a,b]_2-(\rho_2(a)f)b)={\rm i}(\rho_1(a)f)b,
  \end{gather}
  where $[\cdot,\cdot]_1=[\cdot,\cdot]$ and $[\cdot,\cdot]_2={\rm i}[\cdot,\cdot]$.

An immediate consequence of equations \eqref{fail_leib_re} and \eqref{fail_leib_im} is the following.
\begin{Proposition}
Let $(A, [\cdot,\cdot], \rho)$ be a CLA. If either $\mathcal{A}_1=(A_{\R}, [\cdot,\cdot], \rho_1)$ or $\mathcal{A}_2=(A_{\R}, {\rm i}[\cdot,\cdot], -\rho_2)$ are RLAs, then $\rho=0$.
\end{Proposition}

In general, $\mathcal{A}_1$ and $\mathcal{A}_2$ satisfy the Jacobi identity but they are not RLAs, since there is an ``error term'' in the Leibniz identity controlled by the anchor of the other algebroid.

We avoid these error terms by working with $\ker\rho_1$ and $\ker\rho_2$. Also note that the multiplication by ${\rm i}$ is a kind of ``isomorphism'' from $\mathcal{A}_1$ to $\mathcal{A}_2$.

In the same way as with CAVBs, we consider the distribution
\[
A^{\real}=\ker \rho_2=\rho^{-1}(\Delta)
\]
and the real rank of a CLA as
\[\operatorname{\text{real-rank}} A\colon \ m\in M\mapsto \dim_{\R} A^{\real}|_m.\]
We also have the distribution given by $\ker\rho_1$. However, note that
\[\ker\rho_1=j_A(A^{\real})=\rho^{-1}({\rm i}\Delta),\]
 where $j_A$ is the multiplication by ${\rm i}$ in the fibers. Moreover,
\begin{equation}\label{id_1}
[A^{\real}, A^{\real}]\subseteq A^{\real},\qquad[A^{\real}, j_A(A^{\real})]\subseteq j_A(A^{\real})\qquad\text{and}\qquad [j_A(A^{\real}), j_A(A^{\real})]\subseteq A^{\real}.
\end{equation}
 Consider
 \[\rho^{\real}=\rho_1|_{A^{\real}}\qquad \text{and}\qquad \rho^{{\rm im}}=-\rho_2|_{j_A(A^{\real})}.\]

As a consequence, we get the following.
 \begin{Proposition}
 If $(A, [\cdot, \cdot], \rho)$ is a CLA with constant real rank, then $(A^{\real}, [\cdot, \cdot]_{|_{A^{\real}}}, \rho^{\real})$ and
 \smash{$(j_A(A^{\real}), j_A[\cdot, \cdot]_{|_{j_A(A^{\real})}}, \rho^{{\rm im}})$} are isomorphic RLAs, with the isomorphism given by the restric\-tion~of~$j_A$. Moreover, \smash{$\rho^{\real}(A^{\real})=\rho^{{\rm im}}(j_A(A^{\real}))=\Delta$}.
 \end{Proposition}

 {\samepage \begin{Definition}
 When $A^{\real}$ is a RLA, we shall call it as the {\em Lie algebroid of real elements}.
 \end{Definition}

 One of the key feature of the Lie algebroid of real elements is that it allows us to associate a~foliation to any CLA with constant real index.
 \begin{Corollary}
 If a CLA $(A, [\cdot, \cdot], \rho)$ has constant real index, then the set of real elements of the image of $\rho$, $\re\rho(A)$ coincides with the image of $\rho(A^{\real})$. Consequently, $\re\rho(A)$ is an integrable distribution.
 \end{Corollary}

 }

 A RLA is said to be {\em regular} when the image of the anchor map is a subbundle of $TM$. CLAs have three associated distributions: the image of the anchor map, $\Delta$ and $D$, see equation \eqref{real_distributions}. Therefore, a regularity condition for CLAs should take these two additional distributions into consideration.
 \begin{Definition}
 A CLA $(A,[\cdot,\cdot], \rho)$ is {\em regular} if $\rho(A)$ is a regular distribution and we say that it is {\em strongly regular} if $D$ and $\Delta$ are regular distributions.
\end{Definition}
 Note that strongly regular CLAs are regular CLAs with constant real rank, and vice versa.
\begin{Remark}
 Another natural invariant associated to a CLA is the {\em order}, which is a $\mathbb{N}$-function defined as follows:
 \[\order A\colon \ m\in M\mapsto \cork_{\R} D|_{m}.\]
The order was previously introduced in \cite{aguero2022complex} for complex Dirac structures and in \cite{tennyson2021exceptional}, under the name ``class'', for exceptional complex structures.
By equation \eqref{sum_identity}, we have that
\begin{equation}\label{orderrr}
\order A|_m=\rerk A|_m+\dim M-\rk A_{\R},
\end{equation}
and so, constant order is equivalent to constant real rank.
\end{Remark}

In general, the distribution $D$ is not integrable, even for strongly regular CLAs. Non-Levi-flat CR structures provide examples of strongly regular CLAs where the distribution $D$ is not integrable.
\begin{Definition}
Let $(A, [\cdot,\cdot], \rho)$ be a CLA (or an RLA) and let $m\in M$ be a point. The {\em isotropy Lie algebra of $A$ at $m$} is the complex (or real) Lie algebra
\[\g_m (A)=\ker\rho_m.\]
\end{Definition}
When $A$ is regular, we obtain a bundle of complex (or real) Lie algebras given by $\g(A)=\ker \rho$.
 Note that
 \begin{align*}
  \g(A)_{\R}= A^{\real}\cap j_A(A^{\real}) \qquad \text{and}\qquad
\g(A^{\real})=\g(j_A(A^{\real}))=\g(A)_{\R}.
   \end{align*}
\subsection{Minimal complex subalgebroid and the class}
Given a CLA $(A,[\cdot,\cdot],\rho)$, we associate the following smooth real distribution:
\[A_{\min}=A^{\real}+j_A(A^{\real}).\]
Since $A_{\min}$ is $j_A$-invariant, it is a complex distribution. In what follows in this article $A_{\min}$ will be considered as a complex distribution.
 Denote by \smash{$[\cdot,\cdot]_{\min}=[\cdot,\cdot]_{|_{A_{\min}}}$} and \smash{$\rho_{\min}=\rho|_{A_{\min}}$}. An~immediate consequence of equation \eqref{id_1} is the following.
 \begin{Proposition}
If $A_{\min}$ is a regular distribution, then the triple $(A_{\min}, [\cdot,\cdot]_{\min}, \rho_{\min})$ is a~complex Lie subalgebroid of $A$, satisfying $(A_{\min})^{\real}=A^{\real}$.
 \end{Proposition}
Assuming, for example, that $A$ is strongly regular, $A_{\min}$ is a vector bundle. One immediate consequence of the previous proposition is that $A_{\min}$ is the minimal complex Lie subalgebroid of $A$ with the same associated RLA as $A$.
\begin{Definition}
We call $A_{\min}$ the {\em minimal complex subalgebroid of $A$}. We say that $A$ is a~{\em minimal~CLA} whenever $A=A_{\min}$.
\end{Definition}
Given a strongly regular CLA $(A,[\cdot,\cdot],\rho)$, we canonically associate to it two CLAs: $A_{\min}$ and~\smash{$(A^{\real})_{\C}$}. These are related by the following canonical exact sequence of CLAs:
\begin{equation}\label{minre}
\xymatrix{
0\ar[r]&\g(A)\ar[r]^{\Phi}&(A^{\real})_{\C}\ar[r]^{\Psi}&A_{\min}\ar[r]&0,}
\end{equation}
where the maps $\Phi$ and $\Psi$ are Lie algebroid morphisms given by
\begin{gather*}
\Phi\colon \ \g(A)\to (A^{\real})_{\C}, \qquad a\mapsto a-{\rm i}j_A a,\\
\Psi\colon \ (A^{\real})_{\C}\to A_{\min}, \qquad a+{\rm i}b\mapsto a+j_A b.\nonumber
\end{gather*}
As a consequence of the exact sequence \eqref{minre}, we obtain the following.
\begin{Corollary}
   Let $A$ be a CLA with constant real rank such that $A_{\min}$ is regular. Then, $A_{\min}\cong (A^{\real})_{\C}$ if and only if $A$ is an involutive structure $($see Examples~$\ref{examples_0})$.
\end{Corollary}
Note that $A_{\min}=\rho^{-1}(\Delta_{\C})$. Since $A$ is strongly regular, $\g(A)$ is a bundle of complex Lie algebras, and we have the following exact sequence:
\[
\xymatrix{
0\ar[r]&\g(A)\ar[r]& A_{\min}\ar[r]&\Delta_{\C}\ar[r]&0.}
\]
Now we extend the notion of type of complex Dirac structures \cite[Definition 4.7]{aguero2022complex} to CLAs.

\begin{Definition}
  The {\em type} of $A$ at the point $m\in M$ is defined as
  \[\type A|_m= \frac{1}{2}(\dim D|_m-\dim \Delta|_m).\]
\end{Definition}
Note that CLAs with type $0$ are precisely those whose anchor map has a real image ${\rho(A)\!=\!D_{\C}}$, where $D$ is a real distribution. Also note that the type takes values between $0$ and \smash{$\frac{\dim M}{2}$} or~\smash{$\frac{\dim M-1}{2}$}, depending on the parity of $\dim M$. The type has a different meaning in the more general context of CLAs.

\begin{Proposition}
  The following identity holds:
\[
  \type A|_m=\rk_{\C}A-\rk_{\C}A_{\min}|_m.
\]
\end{Proposition}
\begin{proof}
  Applying the rank-nullity theorem on $\rho$ and $\rho_{\min}$, we have
\begin{gather*}
\rk_{\C}\rho(A)+\rk_{\C}\g(A)=\rk_{\C}A\qquad \text{and} \qquad
\rk_{\C}\Delta_{\C}+\rk_{\C}\g(A)=\rk_{\C}A_{\min}.
\end{gather*}

Subtracting the previous expression, we obtain the desired identity.
\end{proof}

\begin{Corollary}\label{type0}
  A CLA has type $0$ if and only if it is minimal.
\end{Corollary}
Thus, the type measures how close is a CLA to being minimal. Furthermore, a strongly regular CLA has constant type, and $\rho(A)$ defines a transverse CR structure on~$M$ (see \cite[Definition~2.74]{aguero}).

We introduce the following invariant.
\begin{Definition}
  The {\em class} of a CLA $A$ at the point $m\in M$ is given by
  \[\class A|_m=\rk \Delta|_m.\]
   We say that a CLA is of {\em class $k$} if it has constant class $k$.
\end{Definition}
\begin{Proposition}We have the following identities:
   \[\class A|_m= \rk_{\C} A_{\min}|_m-\rk_{\C} \g(A)|_m=\rk_{\R} A^{\real}|_m-\rk_{\R} \g(A)_{\R}|_m.\]
\end{Proposition}
While the type measures how close a CLA is to being minimal, the class measures how close it is to being ``real-degenerate'', i.e., when $A^{\real}=\g(A)_{\R}$. Note that the class is bounded below by $0$ but it is not bounded above. We shall see in Examples~\ref{examples1-1} and~\ref{examples1-2-6} that complexifications provide examples of CLAs with arbitrarily high class.
A CLA has class $0$ at a point $m\in M$ if and only if $\Delta|_m=0$. CR structures and holomorphic Lie algebroid are examples of CLAs of class $0$, as we shall see in Examples~\ref{examples1-1} and~\ref{examples1-2-6}. It is easy to see that the only constant real rank CLAs with type $0$ and class $0$ are bundles of complex Lie algebras.

The real rank, the type and the class are related by the following.
\begin{Proposition} The following identity holds:
  \[
  2\type+\class+\rerk =\rk_{\R}A_{\R}.
  \]
  \end{Proposition}
  \begin{proof}
  Applying the rank-nullity theorem on $\rho$ and $\rho_2$ we have the following:
  \begin{align*}
    2\type&=2(\rk_{\R}D-\rk_{\C}\rho(A))\\
       &=2(\rk_{\R}A_{\R}-\rk_{\R}A^{\real})-2(\rk_{\C}A-\rk_{\C} \g(A))\\
       &=\rk_{\R}A_{\R}-2\rerk+\rk_{\R}\g(A)_{\R}\\
       &=\rk_{\R}A_{\R}-\class-\rerk.
\tag*{\qed}
\end{align*}
\renewcommand{\qed}{}
\end{proof}

  \begin{Example}[involutive structures] \label{examples1-1}
 Given an involutive structure $(E,\rho, [\cdot,\cdot])$, we have the following:
\begin{align*}
  E^{\real} =\Delta=E\cap TM,\qquad
  E_{\min} =\Delta_{\C}.
\end{align*}
Moreover,
\[\rerk E=\rk \Delta,\qquad \type E=\rk E-\rk\Delta\qquad\text{and}\qquad\class E=\rk \Delta.\]
 Note that CR structures have trivial real distribution. In fact, they are the only CLAs with this property. Indeed, assume that $A^{\real}=0$, then $A$ is an involutive because $\g(A)=0$ and $A^{\real}=A\cap TM=0$.
 \end{Example}
\begin{Proposition}
  Let $A$ be a CLA. Then, $A^{\real}=0$ if and only if $A$ is a CR structure.
\end{Proposition}

 \begin{Examples}\label{examples1-2-6}\quad
 \begin{enumerate}\itemsep=0pt
\item Complex Dirac structures:
In the case of a complex Dirac structure $L\subseteq \T_{\C}M$, see \cite{aguero2022complex}, we have
 \[L^{\real}=L\cap \bigl(TM\oplus T^*_{\C}M\bigr).\]

Any complex Dirac structure $L$ with constant order has associated a real Dirac structure~\smash{$\widehat{L}$}. By the proof of \cite[Theorem 5.1]{aguero2022complex}, we have the following exact sequence of Lie algebroids:
\[
\xymatrix{
0\ar[r]&L\cap T^* M\ar[r]& L^{\real}\ar[r]^{\Phi}&\widehat{L}\ar[r]&0,}
\]
where
\begin{align*}
\Phi\colon\ L^{\real}\to \widehat{L},\qquad
X+{\rm i}\xi+\eta\mapsto X+\xi.
\end{align*}
In particular, if $L$ is a generalized complex structure, then $L^{\real}\cong \hat{L}$. By equation \eqref{orderrr}, we~have
\[\rerk L=\order L+\dim M.\]
In general, complex Dirac structures are not minimal algebroids. By Corollary \ref{type0}, only complex Dirac structures with type $0$ are minimal.
\item Holomorphic Lie algebroids: A holomorphic Lie algebroid $(A,[\cdot,\cdot], \rho)$ is seen as a CLA via its associated CLA $(A, [\cdot,\cdot]_{1,0},\rho_{1,0})$, see \cite{laurent2008holomorphic}. Note that $\rho_{1,0}(A)\subseteq T_{0,1}M$ and so $\Delta=0$. Hence, holomorphic Lie algebroids have class $0$
\item Bundle of complex Lie algebras: In this case, the anchor map is zero, so $\Delta=0$. Moreover, $A^{\real}= A_{\R}$ and $A_{\min}=A$.
\item Complexification of Lie algebroids: Let $(A, [\cdot, \cdot], \rho)$ be a RLA and let $(A_{\C}, [\cdot, \cdot]_{\C}, \rho_{\C})$ be its complexification.
Note that $\Delta=D=\rho(A)$ and so $A_{\C}$ is minimal. Furthermore,
\[(A_{\C})^{\real}=A\oplus {\rm i}\ker \rho.\]
As a consequence, \begin{align*}
 \rerk A_{\C}=\rk_{\R} A+\rk_{\R}\g(A),\qquad
  \class A_{\C}=\rk_{\R} A-\rk_{\R}\g(A).
\end{align*}

In particular, $A$ has constant real rank if and only if $A$ is a regular RLA.
\item Complex Poisson bivectors: A complex Poisson bivector $\pi\in \Gamma\bigl(\wedge^2 T_{\C}M\bigr)$ decomposes as $\pi=\pi_1+{\rm i}\pi_2$, where $\pi_1$ and $\pi_2$ are real bivectors (not necessarily Poisson). Denote the anchor map of \smash{$\bigl(T_{\C}^*M\bigr)_{\pi}$} by $\rho$. Then, $\rho$ is expressed in terms of $\pi_1$ and $\pi_2$ in the following way:
\[\rho(\xi+{\rm i}\eta)=\pi(\xi+{\rm i}\eta)=\pi_1(\xi)-\pi_2(\eta)+{\rm i}(\pi_2(\xi)+\pi_1(\eta)),\]
where $\xi, \eta\in T^*M$.
So, the associated RLA is given by
\[\bigl(\bigl(T_{\C}^*M\bigr)_{\pi}\bigr)^{\real}=\{ \xi+{\rm i}\eta\mid \pi_2(\xi)+\pi_1(\eta)=0\}\]
with anchor map $\rho^{\real}(\xi+{\rm i}\eta)=\pi_1(\xi)-\pi_2(\eta)$. So we can associate a real foliation to a~complex Poisson bivector, which is not necessarily a symplectic foliation.
The RLA \smash{$\bigl(\bigl(\bigl(T_{\C}^*M\bigr)_{\pi}\bigr)^{\real}, [\cdot,\cdot]^{\real}, \rho^{\real}\bigr)$} and its associated foliation are subject of future work.
\end{enumerate}
  \end{Examples}
Consider two CLAs $(A_1, [\cdot,\cdot]_1, \rho)$ and $(A_2, [\cdot,\cdot]_2, \widetilde{\rho})$ over the same manifold $M$, and a Lie algebroid morphism $\psi\colon A_1\to A_2$. Decompose the anchor maps into their real and imaginary components $\rho=\rho_1+{\rm i}\rho_2$ and \smash{$\widetilde{\rho}=\widetilde{\rho}_1+{\rm i}\widetilde{\rho}_2$}. Note that \smash{$\psi\circ j_{A_1}=j_{A_2}\circ \psi$}, and \smash{$\rho_1=\widetilde{\rho}_1\circ\psi$} and \smash{$\rho_2=\widetilde{\rho}_2\circ\psi$}. As a consequence of these relations, we have the following.
\begin{Proposition}
A morphism $\psi$ between two CLAs $A_1$ and $A_2$, satisfies the following:
\[\psi((A_1)^{\real})\subseteq (A_2)^{\real} \qquad\text{and}\qquad \psi((A_1)_{\min})\subseteq (A_2)_{\min}.\]
In particular, {$\psi|_{(A_1)^{\real}}$} and {$\psi|_{(A_1)_{\min}}$} are Lie algebroid morphisms in both cases. Moreover, both the real rank and the class are invariants under Lie algebroid automorphisms.
\end{Proposition}
\begin{Remark}
   Making a parallel with complex Dirac structures: given a complex Dirac structure $L\subseteq\T_{\C}M$ with constant real index, we associate to it two Lie algebroids: the RLA $K=\re\bigl(L\cap\overline{L}\bigr)$ and the CLA $K_{\C}$. Note that we are strongly using the fact that $L$ is a~subbundle of~$\T_{\C} M$, where a natural complex conjugation is available. Although, in general, a~CLA~$(A,[\cdot,\cdot], \rho)$ is not naturally realized as a subbundle of the complexification of a~Courant or another Lie algebroid, the image of the anchor map $\rho(A)$ is a distribution of $T_{\C}M$. Thus, we pass the information of the real elements of $\rho(A)$ via the inverse image of the anchor map. Therefore, instead of $K$ we consider $A^{\real}=\rho^{-1}(\Delta)$ and instead of $K_{\C}$ we consider~${A_{\min}=\rho^{-1}({\Delta_{\C}})}$.
\end{Remark}

\subsection{Conjugate CLA}
We recall that given a complex vector bundle $A$, its {\em conjugate} bundle $\cA$ is the vector bundle defined by the following modification on the rule of scalar multiplication:
\[z\cdot_{\cA}e=\overline{z}\cdot_{A} e\]
for all $z\in \C$ and $e\in A$, see \cite{milnor1974characteristic}.

Let $(A, \rho, [\cdot,\cdot])$ be a CLA. The map $\overline{\rho}$ becomes $\C$-linear when defined on $\cA$. The bracket $[\cdot,\cdot]\colon \Gamma\bigl(\overline{A}\bigr)\times \Gamma\bigl(\overline{A}\bigr)\to \Gamma\bigl(\overline{A}\bigr)$ is well-defined and remains $\C$-bilinear. Then, we have the following.

\begin{Proposition}
The triple $\bigl(\cA, \overline{\rho}, [\cdot,\cdot]\bigr)$ is a CLA.
\end{Proposition}
\begin{proof}
Since the bracket remains unchanged, the Jacobi identity is satisfied. It remains to verify the Leibniz identity. For $e_1, e_2\in \Gamma\bigl(\cA\bigr)$ and $f\in \CMC$, we have
\begin{align*}
\bigl[e_1, f\cdot_{\cA} e_2\bigr]&{}=\bigl[e_1, \overline{f}\cdot_{A}e_2\bigr]=f[e_1, e_2]+\bigl(\rho(e_1)\overline{f}\bigr)\cdot_{A}e_2\\
&{}=f[e_1, e_2]+ \overline{\rho(e_1)\overline{f}}\cdot_{\cA}e_2=f[e_1, e_2]+\bigl(\overline{\rho(e_1)}f\bigr)\cdot_{\cA}e_2.
\tag*{\qed}
\end{align*}
\renewcommand{\qed}{}
\end{proof}

We refer to the triple $\bigl(\cA, \overline{\rho}, [\cdot,\cdot]\bigr)$ as the {\em conjugate complex Lie algebroid}.
\section{Local structure}\label{locstr}
In this section, we study the local structure of CAVBs and CLAs with constant real rank. Our main tools for providing these local structures come from the techniques developed in \cite{bursztyn2019splitting}. In~Appendix~\ref{preliminarieslocstr}, we recall the properties of normal vector bundles and Euler-like vector fields, which play a key role in this section.
\subsection{Local description of involutive CAVBs}
Involutive RAVBs inherit the structure of an almost Lie algebroid \cite[Proposition 3.17]{bursztyn2019splitting}. For~involutive CAVBs, we have the same result and the proof is an adaptation of the real case to the complex case.
\begin{Proposition}
A real or complex anchored vector bundle $(A,\rho)$ is involutive if and only if there exists a bracket on $\Gamma(A)$ making $A$ into a real or complex almost Lie algebroid, respectively.
\end{Proposition}

Let $(A,\rho)$ be an involutive CAVB with constant real rank. Consider a bracket $[\cdot,\cdot]$ that makes~$(A,\rho)$ into an almost Lie algebroid. Then, we have an operator $D_{\sigma}\colon \Gamma(A)\to \Gamma(A)$, given by $D_{\sigma}\tau=[\sigma,\tau]$. As mentioned in \cite{bischoff2020deformation}, we have the map
\begin{gather*}
\widetilde{\rho}\colon \ \Gamma (A^{\real}) \to \mathfrak{aut}(A,\rho),\qquad
\widetilde{\rho}(\sigma)=(D_{\sigma}, \rho(\sigma))
\end{gather*}
that fits into the following diagram:
\begin{equation}{\label{lifting_aut}}
\xymatrix{
& \mathfrak{aut}(A,\rho)\ar[d]\\
\Gamma(A^{\real})\ar[ru]^{\widetilde{\rho}}\ar[r]^{\rho}& \mathfrak{X}(M).}
\end{equation}
Note that, the map $\widetilde{\rho}$ lifts to a CAVB automorphism but not necessarily to a complex almost Lie algebroid automorphism.
We shall use the following identification of vector bundles.
\begin{Lemma}\label{identification_cx}
Let \smash{$N\xhookrightarrow{\iota}M$} be a submanifold. Then, the vector bundles ${\nu(T_{\C}M,T_{\C}N)\!\xrightarrow{\pr_{T_{\C}N}}\! T_{\C}N}$ and \smash{$T_{\C}\nu(M,N)\xrightarrow{T_{\C}\pr_N} T_{\C}N$} are isomorphic as real vector bundles.
\end{Lemma}
\begin{proof}
  We recall that $T_{\C}M=TM\times_M TM$ and so $ T(T_{\C}M)=TTM\times_{TM} TTM$. Consequently,
  \[ T(T_{\C}M)|_{T_{\C}N}=(TTM\times_{TM} TTM)|_{TN\times_{N}TN}=TTM|_{TN}\times_{TM|_N} TTM|_{TN}.\]
  Note that
  \begin{align*}
  \nu(T_{\C}M,T_{\C}N)&{}=(T(T_{\C}M)|_{T_{\C}N})/T(T_{\C}N)\\
  &{}=\bigl(TTM|_{TN}\times_{TM|_N} TTM|_{TN}\bigr)/TTN\times_{TN} TTN\\
  &{}\cong \nu(TM,TN)\times_{\nu(M,N)}\nu(TM,TN)\\
  &{}\cong T_{\C}\nu(M,N).
\tag*{\qed}
\end{align*}
\renewcommand{\qed}{}
\end{proof}

\begin{Remark}
   Using the identification of Lemma \ref{identification_cx}, we see that if \smash{$N\xhookrightarrow{\iota} M$} is a transversal to $\rho=\rho_1+{\rm i}\rho_2$, then the anchor map is a map of pairs $\rho\colon \bigl(A,\iota^{!}A\bigr)\to (T_{\C}M, T_{\C}N)$, and $\nu(\rho)\colon \nu\bigl(A,\iota^{!}A\bigr)\to \nu(T_{\C}M, T_{\C}N)\cong T_{\C}\nu(M,N)$ is given by
  \[
  \nu(\rho)=\nu(\rho_1)+{\rm i}\nu(\rho_2).
  \]
\end{Remark}
\begin{Lemma}
   If $X$ is an Euler-like vector field along $N$ $($see Definition $\ref{defTN-EL})$, then its complex tangent extension $X_{T_{\C}}$ $($see Definition $\ref{cx_tangent_extension})$ is Euler-like along $T_{\C}N$.
\end{Lemma}

\begin{proof}
   It is easy to see that $X_{T_{\C}}|_{T_{\C}N}=0$. Let $\mathcal{E}$ denote the Euler vector field associated to~${\nu(M,N)\rightarrow N}$. Since the complexified tangent lift of the scalar multiplication is the scalar multiplication in $T_{\C}\nu(M,N)\to T_{\C}N$, we note that $\mathcal{E}_{T_{\C}}$ is an Euler vector field. Using the identification of Lemma \ref{identification_cx}, we see that $\nu(X_{T_{\C}})=\nu(X)_{T_{\C}}=\mathcal{E}_{T_{\C}}$.
\end{proof}

We now apply the techniques of \cite{bursztyn2019splitting} to the case of CAVBs.
\begin{Proposition}\label{splitting_cavb}
Let $(A,\rho)$ be an involutive CAVB with constant real rank and let \smash{$N\xhookrightarrow{\iota} M$} be a submanifold transverse to $\rho^{\real}$. Then
\begin{enumerate}\itemsep=0pt
\item[$1.$] There exists $\widetilde{X}\in \mathfrak{aut}(A,\rho)$ such that $\widetilde{X}|_{\iota^! A}=0$ and $X={(pr_{TM})}_* \widetilde{X}$ is Euler-like along $N$.
\item[$2.$] Each choice of vector field $\widetilde{X}$ as in $(1)$ determines an isomorphism of CAVBs
\[
\widetilde{\psi}\colon\ p^! \iota^! A \to A|_U,
\]
whose base map is the tubular neighborhood embedding $\psi\colon \nu_N\to U\subseteq M$ $($see Definition~$\ref{defTN-EL})$.
\end{enumerate}
\end{Proposition}

\begin{proof}
(1) The proof is similar to \cite[Theorem 3.13]{bursztyn2019splitting}, we just have to work on $T_{\C}M$ instead of~$TM$.

   (2) By Lemma \ref{existence_esp_section}, consider $\epsilon\in \Gamma(A^{\real})$ such that $\rho(\epsilon)=X$ is Euler-like along $N$ and also consider the lifting of equation \eqref{lifting_aut}, $\widetilde{X}\in \aut(A, \rho)$, given by $\widetilde{X}=\widetilde{\rho}(\epsilon)$.
  Note that $\widetilde{X}$ is $\rho$-related to $X_{T_{\C}}$, by the description of $\aut(A,\rho)$.
   Then \smash{$\nu\bigl(\widetilde{X}\bigr)$} is $\nu(\rho)$-related to $\nu(X_{T_{\C}})=\mathcal{E}_{T_{\C}}$, the last one is an Euler vector field and by Lemma \ref{fiberwise_isom}, $\nu(\rho)$ is a fibre-wise isomorphism. Consequently,~\smash{$\widetilde{X}$}~is an~Euler-like vector field along $\iota^! A$.
   Consider $\Phi_t$ and $\widetilde{\Phi}_t$ the flow of $X$ and $\widetilde{X}$, respectively. Since both are Euler-like, the maps $\lambda_t=\Phi_{-\log(t)}$ and \smash{$\widetilde{\lambda_t}=\widetilde{ \Phi}_{-\log(t)}$} are defined at $t=0$. Note that $ \lambda_t\circ \psi=\psi\circ \kappa_t$, in the case of $t=0$ this means that $\lambda_0\circ\psi=p\circ \iota$, where $\psi\colon \nu_N\to U\subseteq M$ is the tubular neighborhood embedding associated to $X$. Recall that the maps $\Tilde{\Phi}_t$ are CLA automorphisms with base $\Phi_t$, so the maps $\Tilde{\lambda}_t$ restricted to $A|_U$ are automorphism of CLA with base $\lambda_t$, for all $t\geq 0$.
   Since $\widetilde{\lambda}_t\circ\widetilde{\lambda}_s=\widetilde{\lambda}_{t+s}$, we have that $\widetilde{\lambda}_t\circ\widetilde{\lambda}_t=\widetilde{\lambda}_{2t}$ and taking limit to $t=0$, we get that $\widetilde{\lambda}^2_0=\widetilde{\lambda}_0$ and so $\widetilde{\lambda}_0$ is a projection. Given that $\Tilde{X}|_{\iota^! A}=0$, we have that $\Tilde{\lambda}_t|_{\iota^! A}={\rm Id}_{\iota^! A}$, so taking limit $\Tilde{\lambda}_0|_{\iota^! A}={\rm Id}_{\iota^! A}$, so $\iota^! A\subseteq \im\bigl(\Tilde{\lambda}_0\bigr)$ and using the fact that
  \begin{equation}\label{eq_inclusion}
    \rho\circ \Tilde{\lambda}_0=T_{\C}\lambda\circ \rho=T_{\C}p\circ T_{\C}\psi^{-1}\circ\rho,
  \end{equation}
  we obtain that $\im\bigl(\widetilde{\lambda}_0\bigr)=\iota^! A$. Consider the map
  \begin{align*} \Psi\colon \ A|_U \to \iota^! A\times T_{\C}\nu_N, \qquad
   a\in A \mapsto \bigl(\widetilde{\lambda}_0(a), T_{\C}\psi^{-1}(\rho(a))\bigr).
   \end{align*}
By equation \eqref{eq_inclusion}, the image of $\Psi$ lies in $p^!\iota^! A$. Finally, note that $\Psi$ is an isomorphism from $A|_U$ to $p^!\iota^! A$, and we define $\widetilde{\psi}=\Psi^{-1}$ as the desired map.
\end{proof}

\begin{Corollary}
Let $(A,\rho)$ be an involutive CAVB and $m\in M$. Assume that $(A,\rho)$ has constant real rank in a neighborhood of $m$. Let \smash{$N\xhookrightarrow{\iota}M$} be a submanifold such that $\Delta|_m\oplus T_m N=T_m M$, and set $P=\Delta{|_m}$. Then, there exists a neighborhood $U$ of $m$ such that $A{|_{U}}$ is isomorphic to~${\iota^! A\times T_{\C}P}$.
\end{Corollary}
\begin{proof}
Since $N$ is completely transversal to $P=\Delta_{|_m}$, we can choose a trivialization of the normal bundle $\nu_N=N\times P$. Then, by Proposition \ref{splitting_cavb}\,(2), we have
$p^!\iota^! A=\iota^! A\times T_{\C}P$.
\end{proof}

\subsection{Local description of CLA}
Consider a CLA $(A, [\cdot,\cdot],\rho)$ with constant real rank. Since $A$ is a CLA, there exists a natural lift
$\widetilde{\rho}\colon \Gamma(A^{\real})\to \mathfrak{aut}_{CLA}(A)$ of the anchor map $\rho$, as in equation~\eqref{lifting_aut}, now defined by the Lie algebroid bracket.
\begin{Theorem}\label{splitting_CLA}
Let $(A,[\cdot,\cdot],\rho)$ be a CLA with constant real rank, and let \smash{$N\xhookrightarrow{\iota}M$} be a submanifold transverse to $A^{\real}$ $($and so transverse to $A)$. Choose a section $\epsilon\in \Gamma(A^{\real})$, such that $\epsilon|_N=0$ and $\rho(\epsilon)\in \mathfrak{X}(M)$ is Euler-like along $N$. Then, there exists a tubular neighborhood embedding $\psi\colon \nu_N\to U\subseteq M$, depending on the choice of $\epsilon$, and a Lie algebroid isomorphism
\[
\widetilde{\psi}\colon\ p^{!}\iota^{!}A\to A|_U
\]
  covering the tubular neighborhood embedding $\psi\colon \nu_N\to U\subseteq M$
\end{Theorem}
\begin{proof}
   The proof follows the same strategy as that of Proposition \ref{splitting_cavb}, now using the lift described above. Thus, we obtain a map $\widetilde{\psi}\colon p^{!}\iota^{!}A\xrightarrow{}A|_U$ that is an isomorphism as CAVBs, we just have to check that $\widetilde{\psi}$ preserves the Lie algebroid bracket. To this end, we use the argument of \cite[Remark 3.20]{bursztyn2019splitting}. Consider the family of maps
   \begin{align*}
   \Psi_t\colon \ A|_U \to \iota^! A\times T_{\C}\nu_N,\qquad
   a\in A \mapsto \bigl(\widetilde{\lambda}_t(a), T_{\C}\psi^{-1}(\rho(a))\bigr).
   \end{align*}
    Note that the image of $\Psi_t$ is $\kappa_t^! \psi^! A$ and consider the maps \smash{$\widehat{\psi}_t=\Psi^{-1}_t$}. Since each map \smash{$\widehat{\psi}_t$} preserves the bracket for all $t>0$, it follows that their limit \smash{$\widehat{\psi}_0=\widetilde{\psi}$} also preserves it. This completes the~proof.
\end{proof}

\begin{Corollary}\label{loc_splitting_CLA}
Let $(A,[\cdot,\cdot],\rho)$ be a CLA and $m\in M$. Assume that $(A,[\cdot,\cdot],\rho)$ has constant real rank in a neighborhood of $m$. Consider a submanifold \smash{$N\xhookrightarrow{\iota}M$} such that $\Delta_{m}\oplus T_m N=T_m M$ and set $P=\Delta_{m}$. Then, there exists a neighborhood $U$ of $m$ such that $A|_{U}$ is isomorphic to~${\iota^! A\times T_{\C}P}$.
\end{Corollary}

Note that $\iota^! A$ has class $0$ at $m$ and $T_{\C}P$ is a minimal CLA (see Examples \ref{examples1-2-6}\,(4)), so the previous description tells us that a CLA with constant real rank splits in a neighborhood of a~point $m$ as a product of a minimal CLA with a CLA of class $0$ at $m$.
The splitting theorem for complex Dirac structures states that a complex Dirac structure with constant order decomposes around a point $m\in M$ as $B$-transformation of the product of a complex Dirac structure having associated Poisson structure vanishing at $m$ (a CLA having class $0$ at $m$) with the complex Dirac structure associated to the presymplectic leaf passing through $m$ (a minimal CLA, see Examples~\ref{examples1-2-6}\,(1)), see~\cite{aguero2022complex}. With this in mind, we can say that minimal CLAs and class $0$ CLA are fundamental pieces for the splitting.

\begin{Remark}
  If \smash{$N\xhookrightarrow{\iota} M$} is a submanifold transversal to $\Delta$, then a normal form $\widetilde{\psi}\colon p^{!}\iota^{!}A\to A|_U$ induces a normal form for $A^{\real}$. Indeed, the~transversality of $N$ with respect to $\Delta$, implies the transversality of $N$ with respect to the anchor map $\rho$. By Lemma~\ref{existence_esp_section}, there exists ${\epsilon\in \Gamma(A^{\real})}$ such that $\rho_1(\epsilon)=X$ is Euler-like along $N$, this is the same as in the proof of~(2). By~Lemma~\ref{restriction_inf_aut}, the restriction of the lift $\widetilde{\rho}$ to $A^{\real}$ produces a lifting $\widehat{\rho}$ for $A^{\real}$. So $\widehat{\rho}(\epsilon)=\widehat{X}$ satisfies the~hypotheses of Theorem~\ref{splitting_ravb}. This determines an isomorphism of RAVBs, $\widehat{\psi}\colon p^!\iota^! A^{\real}\to A^{\real}|_U$ with the same~base map $\psi$ of Proposition~\ref{splitting_cavb}.
\end{Remark}

Assuming the conditions of Corollary \ref{loc_splitting_CLA}, suppose further that $A$ has constant real rank and constant type around a point $m\in M$. Then $\Delta$ is a regular distribution and in a neighborhood $U\subseteq M$ of $m$, we have that
$TN|_{U\cap N}\oplus \Delta|_{U\cap N}=TM|_{U\cap N}$ and so $\iota^! A$ is of class $0$. Hence, we obtain the following.
\begin{Corollary}
If $(A,[\cdot,\cdot],\rho)$ is a strongly regular CLA, $m\in M$ and $P=\Delta_m$, then in a~neighborhood of $m$, the bundle $A$ splits as the product of a CLA of class $0$ with $T_{\C}P$. Moreover, if $(A,[\cdot,\cdot],\rho)$ is a minimal CLA with constant real rank, then around $m$, we have that $A$ is isomorphic to the product of a bundle of complex Lie algebras with $T_{\C}P$.
\end{Corollary}

\begin{Example}[ivolutive structures]
We apply the splitting Theorem \ref{splitting_CLA} for CLAs to an involutive structure $E$.
\end{Example}

\begin{Proposition}\label{loc_str_inv_st}
 Let $(E,\rho, [\cdot,\cdot])$ be an involutive structure and assume that $E^{\real}=\Delta$ is regular. Consider $m\in M$ and a submanifold \smash{$N\xhookrightarrow{\iota}M$} passing through $m$, which is completely transversal to $\Delta$, i.e., $\Delta_m\oplus T_m N=T_m M$. Then, there exists a neighborhood $U$ of $m$ such that
\[
E|_U\cong \iota^! E\times T_{\C}P,
\]
where $P=\Delta_m$ and $\iota^{!} E$ is a CR structure over $U\cap N$.
\end{Proposition}
\begin{proof}
By Theorem \ref{splitting_CLA}, there exists a neighborhood $U$ of $m$ such that $E|_U\cong \iota^! E\times T_{\C}P$. Since $\Delta$ is regular, shrinking $U$ if necessary, we have that $TN|_{U\cap N}\oplus \Delta|_{U\cap N}=TM|_{U\cap N}$. As~a~consequence, $\iota^! E$ has no real elements, i.e., \smash{$\iota^! E\cap \overline{\iota^! E}=0$} and so it is a CR structure.
\end{proof}

 As a consequence, we recover a classical result of Nirenberg about the local description of involutive structures.
\begin{Corollary}[{\cite[Theorem 1]{nirenberg1958complex}}]
   Let $(E,\rho, [\cdot,\cdot])$ be an involutive structure such that $\Delta=E\cap TM$ and $D=\bigl(E+\overline{E}\bigr)\cap TM$ are regular involutive distributions. Then, for any point $m\in M$, there exists a neighborhood $U$ of $m$ and a coordinate system $(x_l, y_l, p_j, q_k)$, with $l=1,\dots, e-r$, $j=1,\dots, r$ and $k=1,\dots, n-2e+r$, where $\dim M=n$, $\rk E=e$ and $\dim \Delta=r$, such that
\[E|_U= {\rm Span}_{\C} \biggl\{\frac{\partial}{\partial \overline{z}_l}, \frac{\partial}{\partial p_j}\biggr\},\]
where $z_l=x_l+{\rm i}y_l$.
\end{Corollary}
\begin{proof}
Choose a submanifold $N\xhookrightarrow{\iota} M$ satisfying that $T_m N\oplus \Delta_m=T_m M$. By Proposition~\ref{loc_str_inv_st}, there exists a neighborhood $U$ of $m$ such that $E|_U\cong \iota^! E\times T_{\C}P$, where $\iota^! E$ is a CR structure. Since $D$ is involutive and $N$ transversal to $D$, we have that \smash{$D'=\re \bigl(\iota^! E + \overline{\iota^! E}\bigr)$} is an involutive distribution on $U\cap N$. Consequently, $\iota^{!}E$ is a Levi-flat CR structure. The generators~\smash{$\frac{\partial}{\partial \overline{z}_l}$} come from the holomorphic foliation associated to $\iota^! E$, while the generators \smash{$\frac{\partial}{\partial p_j}$} come from~$T_{\C}P$.
\end{proof}

\begin{Example}[complex infinitesimal actions]
Let $\g$ be a complex Lie algebra and $\chi\colon \g\to \Gamma(T_{\C}M)$ a complex infinitesimal action. Decompose $\chi$ as $\chi=\chi_1+{\rm i}\chi_2$, where $\chi_1,\chi_2\colon \g \to \Gamma(TM)$ are $\R$-linear. Assume that $L(\g)^{\real}$ is a trivial vector bundle. Then, we have
\[L(\g)^{\real}=\ker \rho_2=M\times \g_{{\rm re}},\]
where $\g_{\rm re}=\ker\chi_2$ is a real Lie algebra contained in $\g$. Moreover, $\chi_1|_{\g_{\rm re}}\to \Gamma(TM)$ is an infinitesimal action on $M$ and so it defines a RLA, that we denote by $L(\g_{\real})$. Let $m\in M$, and let $S$ be an orbit of $L(\g_{\real})$ passing through $m$. Choose a submanifold \smash{$N\xhookrightarrow{\iota}M$} passing through~$m$, which is completely transversal to $S$. Then, by Corollary \ref{loc_splitting_CLA}, there exists a neighborhood $U$ of~$m$ such that
\[L(\g)|_{U}\cong \iota^! L(\g)\times T_{\C}S.\]
In case the map $\chi$ maps into complete vector fields, the infinitesimal action $\chi_1|_{\g_{\rm re}}$ is complete and so by Lie--Palais theorem, this infinitesimal action integrates to a Lie group action $\sigma\colon G_{\rm re}\to \Diff(M)$, where $G_{\rm re}$ is the Lie group integrating $\g_{\rm re}$. In this case, $S$ is also the orbit of the action on the point~$m$.
\end{Example}

Given a RLA $(A, [\cdot, \cdot], \rho)$, its complexification $(A_{\C}, [\cdot, \cdot]_{\C}, \rho_{\C})$ has constant real rank if and only if $A$ is a regular RLA; see Examples \ref{examples1-2-6}\,(4). So Corollary \ref{loc_splitting_CLA} does not directly apply in this case. However, the conclusions of Corollary \ref{loc_splitting_CLA} are still true for the complexification of any RLA and this follows from the local structure of the original $A$. Indeed, consider a~submanifold \smash{$N\xhookrightarrow{\iota}M$} transversal to $\rho_{\C}$ and a point $m\in M$. Then, $N$ is also transversal to $\rho$ and by~\mbox{\cite[Corollary 4.2]{bursztyn2019splitting}}, there exists a neighborhood $U$ of $p$ and an isomorphism of Lie algebroids $\Psi\colon A|_U\to \iota^! A\times T_{\C}P$, where $P=\Delta|_m$. Note that \smash{$\bigl(\iota^! A\bigr)_{\C}=\iota^!(A_{\C})$} and so, the complexification of $\Psi$
\[\Psi_{\C}\colon\ (A_{\C})|_U\to \iota^!(A_{\C})\times T_{\C}P\]
gives the desired local isomorphism. Summarizing, we have the following.

 \begin{Proposition}
Let $(A,[\cdot,\cdot],\rho)$ be a RLA and let $m\in M$. Consider $(A_{\C},[\cdot,\cdot]_{\C},\rho_{\C})$ with the usual CLA structure given by the complexification of $A$ and let $N\xhookrightarrow{\iota}M$ be a submanifold such that $\Delta_{m}\oplus T_m N=T_m M$. Denote by $P=\Delta_{m}$. Then, there exists a neighborhood $U$ of $m$ such that $(A_{\C})|_{U}$ is isomorphic to $\iota^! (A_{\C})\times T_{\C}P$.
\end{Proposition}

\begin{Remark}
Applying Corollary \ref{loc_splitting_CLA} to complex Dirac structures with constant order, we obtain a different result from the one presented in \cite{aguero2022complex}.
The reason is that the group of Lie algebroid automorphisms is larger than the group of automorphisms of complex Dirac structures, so the local description provided here is weaker than the one given in \cite{aguero2022complex}, in terms of automorphisms. However, in \cite{aguero2022complex} we assume that $L$ has constant real index, a condition that is not assumed here.
\end{Remark}
\section{Complex sum of skew-algebroids}\label{cxsum}
Let $A$ be a real vector bundle. Given two brackets $[\cdot,\cdot]_1$ and $[\cdot,\cdot]_2$ on $\Gamma(A)$, we say that they {\em commute} or that they are {\em compatible} if the following identity holds:
\begin{equation}\label{commute_id}
  ([a,[b,c]_1]_2+\mathrm{c.p.})+([a,[b,c]_2]_1+\mathrm{c.p.})=0
\end{equation}
for all $a,b,c\in \Gamma(A)$.
\begin{Definition}
Let $A\to M$ be a real vector bundle and $\mathcal{A}_1=(A,[\cdot,\cdot]_1, \rho_1)$ and $\mathcal{A}_2=(A,[\cdot,\cdot]_2, \rho_2)$ be two structures of skew-algebroids on $A$. The {\em complex sum} of $\mathcal{A}_1$ and $\mathcal{A}_2$ is the CAVB $(A_{\C}, [\cdot,\cdot], \rho)$, with bracket defined by
\[[a,b]=[a,b]_1+{\rm i}[a,b]_2\]
 for $a,b \in \Gamma(A)$ and then extended by $\C$-linearity to $\Gamma(A_{\C})$, and with anchor map
 \[\rho(a+{\rm i}b)=\rho_1(a)-\rho_2(b)+{\rm i}(\rho_1(b)+\rho_2(a)).\]
\end{Definition}
\begin{Proposition}\label{properties_cxsum}
Let $\mathcal{A}_1$ and $\mathcal{A}_2$ be as above. Then
\begin{enumerate}\itemsep=0pt
\item[$1.$] The complex sum always satisfies the Leibniz identity, and so it is a complex skew-algebroid.
   \item[$2.$] The complex sum is an almost CLA $($the anchor map preserves the brackets, see Appendix~$\ref{preliminaries0})$ if and only if
  \[\rho_1[a,b]_2+\rho_2[a,b]_1=[\rho_2(a),\rho_1(b)]+[\rho_1(a), \rho_2(b)].\]
   \item[$3.$] The complex sum is a CLA if and only if the skew-algebroids commute and have the same Jacobiators ${\rm Jac}[\cdot,\cdot]_1={\rm Jac}[\cdot,\cdot]_2$.
\end{enumerate}
\end{Proposition}
\begin{proof}
(1) and (2): It is enough to prove this for real elements.

(3): The commutativity of the brackets and the equality of the Jacobiator of $[\cdot,\cdot]_1$ and $[\cdot,\cdot]_2$ are together equivalent to the Jacobi identity for the bracket $[\cdot,\cdot]$.
\end{proof}

As a consequence of the previous proposition, we have the following.

\begin{Definition}
 A pair of skew-algebroids structures $(\mathcal{A}_1,\mathcal{A}_2)$ over the same vector bundle $A$ is called a {\em complex matched pair} if they satisfy the hypothesis of Item (3) of Proposition \ref{properties_cxsum}.
\end{Definition}
\begin{Remark}
   There is a definition of matched pairs for real Lie algebroids, which can be extended to CLAs. We shall refer to this extension as {\em $\C$-extended matched pairs}. Our definition is slightly different. While a $\C$-extended matched pairs consists of a pair of CLAs $(\mathcal{A},\mathcal{B})$ over the same base manifold $M$, such that $\mathcal{C}=\mathcal{A}\oplus \mathcal{B}$ is a CLA over $M$ and $\mathcal{A}$ and $\mathcal{B}$ are subalgebroids of~$\mathcal{C}$, see \cite[Definition~4.1]{mokri1997matched}, our definition could be seen as a pair $(\mathcal{A}_1,\mathcal{A}_2)$ of skew-algebroids over the same vector bundle $A\to M$ such that ``$\mathcal{A}_1+{\rm i}\mathcal{A}_2$'' is a CLA, as defined above. Moreover,~$\mathcal{A}_1$ and~$\mathcal{A}_2$ are not necessarily subalgebroids.
\end{Remark}
\begin{Examples}\quad
\begin{enumerate}\itemsep=0pt
\item The complexification of a RLA $(A, [\cdot, \cdot], \rho)$ is the complex sum of $(A, [\cdot, \cdot], \rho)$ with \linebreak ${(A, [\cdot, \cdot]=0, 0)}$.
\item The complex double of a RLA $(A, [\cdot, \cdot], \rho)$ is the complex sum of $(A, [\cdot, \cdot], \rho)$ with itself. In~this~case
\[A^{\real}=\{a+{\rm i}(\gamma-a)\mid a\in A, \, \gamma\in \ker\rho\}.\]

\end{enumerate}
\end{Examples}

\begin{Theorem}\label{cx_uniqueness}
Let $A$ be a real vector bundle over $M$. Any structure of CLA on $A_{\C}$ comes from a complex sum of a complex matched pair of skew-algebroids.
\end{Theorem}
\begin{proof}
Consider the $\R$-linear maps $\Tilde{\rho}_1,\Tilde{\rho}_2\colon A_{\C}\to TM$ given by the decomposition $\rho=\Tilde{\rho_1}+{\rm i}\Tilde{\rho_2}$. Note~that
\[\rho(a+{\rm i}b)=(\Tilde{\rho}_1(a)-\Tilde{\rho}_2(b))+{\rm i}(\Tilde{\rho}_1(b)+\Tilde{\rho}_2(a)).\]
Now consider the $\R$-linear maps $\rho_1=\Tilde{\rho}_1|_{A}$ and $\rho_2=\Tilde{\rho}_2|_{A}$. Then
\[\rho(a)=\rho_1(a)+{\rm i}\rho_2(a) \qquad \forall\ a\in A.\]

Now consider the $\R$-linear maps $B_1, B_2\colon \Gamma(A_{\C})\times \Gamma(A_{\C})\to \Gamma (A)$ defined as follow:
\[[a_1+{\rm i} a_2, b_1+{\rm i}b_2]=B_1(a_1+{\rm i} a_2, b_1+{\rm i}b_2)+{\rm i}B_2(a_1+{\rm i} a_2, b_1+{\rm i}b_2).\]
In the same way, as with the decomposition of $\rho$, we get that
\begin{equation}\label{def_brackets}
  [a_1,a_2]=B_1(a_1, a_2)+{\rm i}B_2(a_1,a_2) \qquad\forall\ a_1,a_2\in \Gamma(A).
\end{equation}

Consider the brackets $[\cdot,\cdot]_1=B_1|_{\Gamma(A)\times\Gamma(A)}$ and $[\cdot,\cdot]_2=B_2|_{\Gamma(A)\times \Gamma(A)}$ on $\Gamma(A)$. The conditions of item (3) of Proposition \ref{properties_cxsum}, follow from examining the real and imaginary parts appearing in the Jacobi identity of $[\cdot,\cdot]$ after decomposing it in terms of $B_1$ and $B_2$. The real part yields $\jc[\cdot,\cdot]_1=\jc[\cdot,\cdot]_2$, while the imaginary part is equivalent to the compatibility condition for the~brackets.

Using equation \eqref{def_brackets} and the Leibniz identity, we have that
\begin{align*}
  B_1(a_1, fa_2)+{\rm i}B_2(a_1, fa_2)&= [a_1, fa_2]
   =f[a_1,a_2]+(\rho(a_1)f)a_2\\
&=fB_1(a_1,a_2)+(\rho_1(a_1)f)a_2+{\rm i}(fB_2(a_1,a_2)+(\rho_2(a_1)f) a_2),
\end{align*}
so both $[\cdot,\cdot]_1$ and $[\cdot,\cdot]_2$ satisfy the Leibniz identity.
The triples $(A,[\cdot,\cdot]_1, \rho_1)$ and $(A,[\cdot,\cdot]_1, \rho_1)$ do not necessarily satisfy the Jacobi identity, so they are skew-algebroids and their complex sum is~$(A_{\C},[\cdot,\cdot],\rho)$ by construction.
\end{proof}

The distributions associated to a complex sum $A_{\C}$ are the following:
\begin{align*}
  &D=(\rho_1(A)+\rho_2(A))_{\C},\\
  &\Delta=\{\rho_1(a)-\rho_2(b)\mid \rho_2(a)+\rho_1(b)=0\ \text{and}\ a,b\in A \}\qquad\text{and}\\
  &(A_{\C})^{\real}=\{a+{\rm i}b\in A_{\C} \mid \rho_2(a)+\rho_1(b)=0\}.
\end{align*}

Let $\mathcal{A}_1$ and $\mathcal{A}_2$ be two skew-algebroid structures on $A$. Then, they have associated fiber-wise linear bivectors $\pi_1$ and $\pi_2$ on $A^*$, respectively (see Appendix \ref{preliminaries0}).
\begin{Proposition}\label{cxlinbivcxmp}
The compatibility condition of the brackets is equivalent to the compatibility condition of the bivectors $[\pi_1,\pi_2]=0$. The condition ${\rm Jac}[\cdot,\cdot]_1={\rm Jac}[\cdot,\cdot]_2$ is equivalent to the condition $[\pi_1,\pi_1]=[\pi_2,\pi_2]$.
\end{Proposition}
\begin{proof}
First, we recall the identity
 \[[\pi_1,\pi_2](f,g,h)=\pi_1(\pi_2(f,g), h)+\pi_2(\pi_1(f,g), h)+\mathrm{c.p.},\]
where $\pi_1,\pi_2\in \mathfrak{X}^2(M)$, $f,g,h\in C^{\infty}(M,\R)$, and we see the multivectors as multi-derivations. The result follows by evaluating this identity on linear and basic functions of $C^{\infty}(A^*,\R)$.
\end{proof}

These equivalences provide geometric interpretations of the conditions defining a complex matched pair.
\subsection{Complex vector fields}
As an application, we study the structure of CLAs associated to complex vector fields, see Example~\ref{examples_0}\,(8).
The following Proposition is a well-known fact, but we provide a proof for completeness:
\begin{Proposition}
  The triple $(A_Z, [\cdot,\cdot], \rho)$ is a CLA.
\end{Proposition}
 \begin{proof}
   It is straightforward to verify that the conditions for being an almost CLA are satisfied. For the second part, note that $A_Z$, as a bundle, is the complexification of $M\times \R$ and that~$(A_Z, [\cdot,\cdot], \rho)$ is the complex sum of the Lie algebroids $A_{X_1}$ and $A_{X_2}$,
   where $Z=X_1+{\rm i}X_2$, $X_1,X_2\in \X(M)$. So, by Proposition \ref{cx_uniqueness}, $A_Z$ is a CLA if and only if the brackets of~$A_{X_1}$ and~$A_{X_2}$ commute, since the equality of the Jacobiators is already satisfied (both have vanishing Jacobiators). The brackets of~$A_{X_1}$ and~$A_{X_2}$ commute if and only if equation \eqref{commute_id} is satisfied
   \[([f,[g,h]_1]_2+c.p)+([f,[g,h]_2]_1+\mathrm{c.p.})=0.\]
   Expanding the terms, we have that
   \begin{align*}
   [f,[g,h]_1]_2+[f,[g,h]_2]_1&{}=
fg(L_{X_1}L_{X_2}h+L_{X_2}L_{X_1}h)-fh(L_{X_1}L_{X_2}g+L_{X_2}L_{X_1}g)\\
&\quad{}-g(L_{X_1}h L_{X_2}f+L_{X_2}h L_{X_1}f)+h(L_{X_1}g L_{X_2}f+L_{X_2}g L_{X_1}f).
     \end{align*}
 The remaining terms are permutations of this formula. Note that all these terms cancel each other out and so the proposition holds.
\end{proof}

Now we study the distribution $(A_Z)^{\real}$. First, note that
\[\rho|_m(z_1+{\rm i}z_2)=z_1 X_1|_m - z_2 X_2|_m + {\rm i}(z_2 X_1|_m + z_1 X_2|_m).\]
By analyzing the different cases, we calculate $(A_Z)^{\real}=\ker\rho_2$. In summary, we obtain the following description:
\[
(A_Z)^{\real}|_m=
\begin{cases}
\C_m, & X_1|_m=X_2|_m=0,\\
(\R\oplus {\rm i}0)|_m, & X_1|_m\neq 0, \ X_2|_m=0,\\
(0\oplus {\rm i}\R)|_m, & X_1|_m=0, \ X_2|_m\neq 0,\\
0_m, & \begin{array}{@{}l@{}} X_1|_m\neq 0, \ X_2|_m\neq 0,\\
    X_1|_m, \  X_2|_m\ \text{are linearly independent},\end{array}\\
\{-ct+{\rm i}t\mid t\in \R\}_m,& X_1|_m\neq 0, \ X_2|_m\neq0\ \text{and}\ X_2|_m=cX_1|_m
\end{cases}
\]
and
\[
(A_Z)_{\min}|_m=
\begin{cases}
0, & X_1|_m\neq 0,\  X_2|_m\neq 0,\  X_1|_m, \ X_2|_m\ \text{are linearly independent},\\
\C_m,& \text{elsewhere}.
\end{cases}
\]
Finally, we have the following.
\begin{Proposition}\label{cxvf}
The distribution $A_Z$ has constant real rank if and only if either
\begin{enumerate}\itemsep=0pt
  \item[$1.$] $X_1=X_2=0$.
  \item[$2.$] $X_1|_m\neq 0$, $X_2|_m\neq 0$ and $X_1|_m$ and $X_2|_m$ are linearly independent.
   \item[$3.$] $X_1$ and $X_2$ never vanish simultaneously, and there exist smooth functions $c\in C^{\infty}(M\setminus {\rm Sing}(X_1))$ and $d\in C^{\infty}(M\setminus {\rm Sing}(X_2))$ such that $c X_1=X_2$ and $d X_2=X_1$.
\end{enumerate}
\end{Proposition}
As a consequence of Proposition \ref{cxvf}, in cases (1) and (2) the foliation associated to $A_Z$ consists of isolated points, while in case (3) by the integral curves of $X_1$ or $X_2$, depending on which one does not vanish.

\begin{Proposition}
 Any structure of CLA over $M\times \C\to M$ comes from a complex vector field and is a complex sum of the RLAs defined by the real and imaginary components of this complex vector field.
\end{Proposition}
\begin{proof}
 Let $(M\times\C, [\cdot,\cdot], \rho)$ be any structure of CLA over $M\times \C$. Consider the complex vector field $Z$ given by $Z_m=\rho(m,1)$ or equivalently $Z=\rho(1)$, where $1\in C^{\infty}(M,\C)$ is the constant function. So, for any $f\in C^{\infty}(M,\C)$, we have $\rho(f)=fZ$. By the Leibniz identity, we have
 \[[1,f]=[1,f\cdot 1]=L_Z f\]
 and thus
 \[[f,g]=[f,g\cdot 1]=g[f,1]+(\rho(f)g)1=fL_Z g-gL_Z f.\]
 Consider $X_1,X_2\in \mathfrak{X}(M)$ such that $Z=X_1+{\rm i}X_2$. By evaluating the anchor map and the bracket on real functions, it follows that
 \[\rho(f)=fX_1+{\rm i}fX_2\qquad \text{and}\qquad [f,g]=[f,g]_1+{\rm i}[f,g]_2,\]
 where $[\cdot,\cdot]_1$ and $[\cdot,\cdot]_2$ are the brackets of the RLAs defined by $X_1$ and $X_2$.
\end{proof}

\subsection{Complex bivectors}
Consider a real bivector $\gamma\in \Gamma\bigl(\wedge^2 TM\bigr)$ and the bracket
\[[\cdot,\cdot]_{\gamma}\colon\ \Gamma(T^*M)\times\Gamma(T^*M)\to \Gamma(T^*M)\]
defined by
\begin{equation}\label{bivector_bracket}
[\alpha,\beta]_{\gamma}=L_{\gamma(\alpha)}\beta-L_{\gamma(\beta)}\alpha-d\gamma(\alpha,\beta)
\end{equation}
   for $\alpha, \beta\in \Gamma(T^*M)$. For a complex bivector $\pi\in \Gamma\bigl(\wedge^2 T_{\C}M\bigr)$, the same formula as in \eqref{bivector_bracket} defines a bracket $[\cdot,\cdot]_{\pi}$ on $\Gamma(T_{\C}M)$. A straightforward verification yields the following.
\begin{Proposition}
  If $\gamma$ is a real or complex bivector, then the triple $(T^*M,[\cdot,\cdot]_{\gamma}, \gamma)$ or $\bigl(T_{\C}^*M,\allowbreak [\cdot,\cdot]_{\gamma}, \gamma\bigr)$ is a skew-algebroid or a complex skew-algebroid, respectively. We denote these structures by \smash{$(T^*M)_{\gamma}$} or~\smash{$\bigl(T_{\C}^*M\bigr)_{\gamma}$}, respectively.
\end{Proposition}

 Consider a complex bivector $\pi=\pi_1+{\rm i}\pi_2\in \Gamma\bigl(\wedge^2 T_{\C}M\bigr)$, where $\pi_1$ and $\pi_2$ are real bivectors. Denote by $\rho$ the anchor map of \smash{$\bigl(T_{\C}^*M\bigr)_{\pi}$}. Then, $\rho$ decomposes as
\[\rho(\xi+{\rm i}\eta)=\pi(\xi+{\rm i}\eta)=\pi_1(\xi)-\pi_2(\eta)+{\rm i}(\pi_2(\xi)+\pi_1(\eta))\]
for $\xi, \eta\in T^*M$.
Moreover, the bracket evaluated on real elements satisfies
\[[\alpha,\beta]_{\pi}=[\alpha,\beta]_{\pi_1}+{\rm i}[\alpha,\beta]_{\pi_2}.\]

\begin{Proposition}\label{cxbivcxmp}
 Let $\pi=\pi_1+{\rm i}\pi_2$ be a complex bivector. Then, $(T_{\C}^*M)_{\pi}$ is the complex sum of $(T^*M)_{\pi_1}$ and $(T^*M)_{\pi_2}$. Moreover, $\pi$ is a complex Poisson bivector if and only if $(T^*M)_{\pi_1}$ and $(T^*M)_{\pi_2}$ form a complex matched pair.
\end{Proposition}
\begin{proof}
   The first part of the proposition follows from what was exposed above. The bivector $\pi$ is Poisson if and only if \smash{$\bigl(T_{\C}^*M\bigr)_{\pi}$} is a CLA. So by item (3) of Proposition \ref{properties_cxsum}, this holds if and only if $(T^*M)_{\pi_1}$ and $(T^*M)_{\pi_2}$ form a complex matched pair.
\end{proof}

\begin{Remark}
In particular, we note that a pair of Poisson bivectors $(\pi_1,\pi_2)$ is a bi-Hamiltonian structure if and only if $(T^*M)_{\pi_1}$ and $(T^*M)_{\pi_2}$ form a complex matched pair. By combining Propositions~\ref{cxlinbivcxmp} and~\ref{cxbivcxmp}, we note that the notion of complex matched pair on a~pair of Lie algebroids $(\mathcal{A}_1,\mathcal{A}_2)$ actually extends the notion of bi-Hamiltonian structures to the context of Lie algebroids. This extension leads to the notion of ``bi-Hamiltonian Lie algebroids''.
\end{Remark}

\appendix

\section{Anchored vector bundles and Lie algebroids}\label{preliminaries0}
In this appendix, we recall some properties of real anchored vector bundles and real Lie algebroids.
\begin{Definition}
   A {\em real anchored vector bundle $($RAVB$)$} is a pair $(A, \rho)$, where $A$ is a vector bundle over a manifold $M$ and $\rho\colon A\to TM$ is a $\R$-bundle morphism called the {\em anchor map}. We~say that a RAVB is {\em involutive} whenever $\rho(\Gamma(A))$ is a Lie subalgebra of $\Gamma(TM)$.
\end{Definition}
Given two RAVBs $(A,\rho_A)$ and $(B,\rho_B)$ over $M$ and $N$, respectively. A {\em morphism of RAVBs} is given by a bundle map $\Phi\colon A\to B$ and a map $\varphi\colon M\to N$ such that the following diagram commute:
\[\xymatrix{A\ar[d]^{\rho_A}\ar[r]^{\Phi} & B\ar[d]^{\rho_B}\\ TM\ar[r]^{T\varphi} & TN.
}\]
We recall that the Lie algebra $\mathfrak{aut}(E)$ of infinitesimal automorphism of a real vector bundle~$E$ are vector fields \smash{$\widetilde{X}\in \X(E)$} whose flow is given by automorphism of $E$. Equivalently, they are given by derivations, that is, pairs $(L,X)$, where $L\colon \Gamma(E)\to \Gamma(E)$ is a $\R$-linear operator and~${X\in \X(M)}$ satisfying $L(f\sigma)=fL(\sigma)+X(f)\sigma$. The relationship is given by the following: given~${\widetilde{X}\in \aut(E)}$, consider the pair $(L,X)$ defined by
\[
X=\pr_* \widetilde{X},\qquad \pr^*\langle \eta, X \rangle=\bigl\langle \pr^*\eta,\widetilde{X} \bigr\rangle\qquad \text{and}\qquad
L(\sigma)^{\uparrow}=\bigl[\sigma^{\uparrow},\widetilde{X}\bigr],
\]
where $\sigma\in \Gamma(E)$, $\eta\in \Omega^{1}(M)$ and $\sigma^{\uparrow}$ is the vertical lift of $\sigma$:
\begin{align*}
   \sigma\in \Gamma(E)\to \sigma^{\uparrow}\in\X(E),\qquad
 \sigma^{\uparrow}(e_m)=\frac{{\rm d}}{{\rm d}t}\bigg\rvert_{t=0}(e_m+t\sigma(m))\qquad\forall\ e_m\in E_m.
 \end{align*}
 Note that this map is injective, since $\sigma^{\uparrow}=0$, implies that $\sigma=0$ just by evaluating in the zero section of $E$.

 The {\em space of automorphisms of a RAVB $(A,\rho)$} is given by
 \[\Aut(A,\rho)=\{(\Phi, \varphi)\in \Aut(A)\mid \rho\circ\Phi=T\varphi\circ \rho\}\]
 and the {\em space of infinitesimal automorphisms of a RAVB $(A,\rho)$} is given by
 \[\mathfrak{aut}(A,\rho)=\big\{ \widetilde{X}\in \aut(A)\mid \widetilde{X}\sim_{\rho} X_T \big\},\]
 where $X_T\in \X(TM)$ is the vector field constructed with the differential of the flow of $X$. Equivalently, the infinitesimal automorphisms of $(A,\rho)$ are given by derivations $(L,X)$ satisfying~that
 \[\rho(L(\sigma))=[X, \rho(\sigma)].\]

The pullback of a RAVB $(A, \rho) $ over $M$ along the map $\varphi\colon N\to M$ is defined as the fibered product:
\[
\xymatrix{\varphi^! A\ar[d]\ar[r] & A \ar[d]^{\rho}\\
TN\ar[r]^{T\varphi} & TM.}
\]
The local description of an involutive RAVB is the following.
\begin{Theorem}[{\cite[Corollary 3.14]{bursztyn2019splitting}}]\label{splitting_ravb}
Let $(A,\rho)$ be an involutive RAVB, $m\in M$, and \smash{$N\xhookrightarrow{\iota} M$} a~submanifold containing $m$ satisfying that $\rho(A|_m)\oplus T_m N=T_m M$. Let $P=\rho(A|_m)$. Then, there exists a neighborhood $U$ of $m$ such that $A|_{U}$ is isomorphic to $\iota^! A\times TP$.
\end{Theorem}
\begin{Definition}A {\em real Lie algebroid $($RLA$)$} over a manifold $M$ is a triple $(A, [\cdot, \cdot], \rho)$, where $(A,\rho)$ is a RAVB over $M$ together with a Lie bracket on $\Gamma(A)$:
\[[\cdot,\cdot]\colon\ \Gamma(A)\times \Gamma(A)\to\Gamma(A)\]
satisfying the Leibniz property
\[[\alpha,f\beta]=f[\alpha,\beta]+\rho(\alpha)(f)\beta\]
 for all $\alpha,\beta\in \Gamma(A)$ and $f\in C^{\infty}(M)$. If the triple $(A,[\cdot,\cdot],\rho)$ satisfies all the previous conditions except the Jacobi identity, it is called a {\em skew-algebroid}, see \cite[Definition 4.1]{grabowski2022nonassociative}; they were originally introduced in \cite{kosmann1990poisson} as differential pre-Lie algebras.
\end{Definition}

Skew-algebroids $A$ are equivalent to fiber-wise linear bivectors on $A^*$, in the following way:
\begin{align*}
\{f\circ p_{A^*},g\circ p_{A^*}\}=0,\qquad
\{l_{\alpha}, f\circ p_{A^*}\}=\rho(\alpha)f,\qquad
\{l_{\alpha},l_{\beta}\}=l_{[\alpha,\beta]},
\end{align*}
where $p_{A^*}\colon A^* \to M$ is the usual bundle projection, $f, g\in C^{\infty}(M,\R)$, $\alpha,\beta\in \Gamma(A)$ and ${l_{\alpha}, l_{\beta}\in C^{\infty}(A^*,\R)}$ are defined as $l_{\alpha}(\tau)=\langle\tau,\alpha\ra$. In case these bivectors are Poisson bivectors, we obtain Lie algebroid structures.

The pullback of RLAs is defined in the same way as the pullback of RAVBs. We recall the local description of RLAs.

\begin{Theorem}[{\cite{bursztyn2019splitting, fernandes2002lie, Dufour2001NormalFF}}]
 Let $(A,[\cdot,\cdot],\rho)$ be a RLA and $m\in M$. Consider a submanifold \smash{$N\xhookrightarrow{\iota}M$} such that $\rho(A|_{m})\oplus T_m N=T_m M$ and denote by $P=\rho(A|_{m})$. Then, there exists a~neighborhood $U$ of $m$ such that $A|_{U}$ is isomorphic to $\iota^! A\times TP$.
\end{Theorem}

The following table summarizes all the algebroids used in this article:
\begin{center}
    \begin{tabular}{@{}c|c@{}}
      Real & Complex\\
      \hline
Real anchored vector bundle (RAVB): & Complex anchored vector bundle (CAVB):\\
$\R$-vector bundle and $\R$-anchor map & $\C$-vector bundle and $\C$-anchor map\\
\hline
Involutive RAVB: & Involutive CAVB:\\
A RAVB where the image of the anchor map & A CAVB where the image of the anchor map\\
 is involutive in $\Gamma(TM)$ & is involutive in $\Gamma(T_{\C}M)$ \\
\hline
 Real skew-algebroid: & Complex skew-algebroid: \\
$\R$-vector bundle, $\R$-anchor map and &
$\C$-vector bundle, $\C$-anchor map and\\
$\R$-bilinear bracket satisfying only Leibniz & $\C$-bilinear bracket satisfying only Leibniz \\
\hline
Real almost Lie algebroid: & Complex almost Lie algebroid \\
It is a real skew-algebroid where &
It is a complex skew-algebroid where\\
 the anchor map preserves the bracket & the anchor map preserves the bracket \\
\hline
Real Lie algebroid (RLA): & Complex Lie algebroid (CLA):\\
It is a real skew-algebroid where &
It is a complex skew-algebroid where\\
 the bracket satisfies Jacobi & the bracket satisfies Jacobi \\
\hline
    \end{tabular}
    \end{center}

\section{Euler-like vector fields and normal forms}\label{preliminarieslocstr}
 In this appendix we recall the properties of normal vector bundles and Euler-like vector fields, that are mainly used in Section~\ref{locstr}.
The {\em Euler vector field} associated to a vector bundle $E\to M$, usually denoted by $\mathcal{E}\in \mathfrak{X}(E)$, is the vector field generated by the one-parameter group $s\mapsto \kappa_{e^{-s}}$, where the maps $\kappa_{t}\colon E\to E$ are given by $\kappa_{t}e=t\cdot e$, $t\in \R$ (in case $t=0$, $\kappa_0$ is the projection to the zero section).

Given a submanifold $N$ of $M$, we denote the normal bundle of $N$ by
\[\nu(M,N)=TM|_{N}/TN,
\]
 when the context is clear we denote $\nu(M,N)$ by $\nu_N$ and the projection map by $p\colon \nu(M,N)\to N$.

Any map of pairs $\varphi\colon (M',N')\to (M,N)$ has an associated map $\nu(\varphi)\colon \nu(M',N')\to \nu(M,N).$
The following lemma is a known fact.
\begin{Lemma}\label{fiberwise_isom}
Let $\varphi\colon M'\to M$ be a smooth map. If \smash{$N\xhookrightarrow{\iota}M$} is a submanifold transverse to $\varphi$ and $N'=\varphi^{-1}(N)$, then $\nu(\varphi)$ is a fiberwise isomorphism.
\end{Lemma}
Consider a vector field $X\in \mathfrak{X}(M)$ tangent to $N$, that is $X\colon (M,N)\to (TM,TN)$. Using the fact that $T\nu_N=\nu_{TN}$ (see \cite[Appendix A]{bursztyn2019splitting}), we have that $\nu(X)$ is a vector field on~$\nu(M,N)$.
\begin{Definition}\label{defTN-EL}
Let $N$ be a submanifold of $M$. A {\em tubular neighborhood embedding} is an embedding
$\psi\colon \nu_N\to M$ that takes $N$ to $N$ and $\nu(\psi)={\rm Id}$, where $\nu(\psi)$ is induced by $\psi\colon (\nu_N,N)\to (M,N)$ (here we are using the identification $\nu(\nu_N,N)=\nu_N$).
A vector field $X\in\mathfrak{X}(M)$ is called {\em Euler-like} (along $N$) if it is complete, $X|_{N}=0$ and $\nu(X)$ is the Euler vector field of $\nu(M,N)$.
\end{Definition}
Given a submanifold $N$ of $M$, there is a one-to-one correspondence between Euler-like vector and tubular neighborhood embedding, see, for example, \cite{bursztyn2019splitting}.
\begin{Lemma}[{\cite[Lemma 3.9]{bursztyn2019splitting}}]\label{existence_esp_section}
Let $(A, \rho)$ be RAVB over $M$, and $N\subseteq M$ a submanifolds transversal to $\rho$.
Then there exists a section $\epsilon\in \Gamma(A)$ with $\epsilon|_N = 0$, such that $X = \rho(\epsilon)$ is Euler-like along~$N$.
\end{Lemma}

Given a vector bundle \smash{$E\xrightarrow{p} M$}, we know that
\[
TE\xrightarrow{Tp} TM \qquad \text{and}\qquad TE\times_E TE\xrightarrow{Tp\times Tp} TM \times _M TM
\]
 are vector bundles. We give the structure of complex vector bundle to the last one and we call it {\em the complexified tangent bundle} \smash{$T_{\C}E\xrightarrow{T_{\C}p} T_{\C}M$}, with the fiber sum and fiber multiplication by scalars given by
\begin{align*}
&(V+{\rm i} W)+_{T_{\C}p} \bigl(V'+{\rm i} W'\bigr)=\bigl(V+_{Tp} V'\bigr)+{\rm i}\bigl(W+_{Tp} W'\bigr),\\
&(\lambda_1+{\rm i}\lambda_2)\cdot_{T_{\C}p} (V+{\rm i}W)=(\lambda_1\cdot_{Tp} V-\lambda_2\cdot_{Tp} W) + {\rm i}(\lambda_2\cdot_{Tp} V + \lambda_1\cdot_{Tp} W),
\end{align*}
where $V$, $V'$, $W$ and $W'\in T_{\C}E$, with $p(V)=p(V')$, $p(W)=p(W')$ and $\lambda_1,\lambda_2\in \R$, and~$+_{Tp}$ and~$\cdot_{Tp}$ denotes the operations of fiber sum and multiplication by scalar in the tangent vector bundle \smash{$TE\xrightarrow{Tp}TM$}. Note that the previous operations make \smash{$T_{\C}E\xrightarrow{T_{\C}p} T_{\C}M$} a complex vector~bundle,
\[
\begin{tikzcd}
\Gamma\bigl(\wedge^k T^*_{\C}M\bigr)\arrow[d, "\widetilde{d}"] \arrow[r, "\widehat{}"]
& \Gamma\bigl(\wedge^k T^*M\bigr)_{\C}\arrow[d, "d_{\C}"] \\
\Gamma\bigl(\wedge^{k+1} T^*_{\C}M\bigr)\arrow[r, "\widehat{}"]
& \Gamma\bigl(\wedge^{k+1} T^*M\bigr)_{\C}.
\end{tikzcd}
\]

\subsection*{Acknowledgements}
The author acknowledges Henrique Bursztyn for proposing the topic and his many valuable observations. The author is also grateful to Hudson Lima, Roberto Rubio and Pedro Frejlich for many fruitful and helpful conversations. We are also thankful to the anonymous referees for their valuable comments and suggestions.

\pdfbookmark[1]{References}{ref}
\LastPageEnding

\end{document}